\documentstyle{amsart}
\input{bull-art}
\newtheorem{th1}{Theorem}
\newtheorem{th1a}[th1]{Theorem}
\newtheorem{th2}[th1]{Theorem}

\newtheorem{th5}[th1]{Theorem}
\newtheorem{th6}[th1]{Theorem}
\newtheorem{th7}[th1]{Theorem}
\newtheorem{rep}[th1]{Theorem}
\newtheorem{perf}[th1]{Theorem}
\newtheorem{classification}[th1]{Theorem}
\newtheorem{sullivan}[th1]{Theorem}
\newtheorem{smale}[th1]{Theorem}
\newtheorem{sing}[th1]{Theorem}
\newtheorem{nowd}[th1]{Theorem}
\newtheorem{asympt}[th1]{Theorem}
\newtheorem{bakdom}[th1]{Theorem}
\newtheorem{ber92}[th1]{Theorem}
\newtheorem{ber92b}[th1]{Theorem}
\newtheorem{ere89}[th1]{Theorem}
\newtheorem{ere90}[th1]{Theorem}
\newtheorem{threlax}[th1]{Theorem}
\newtheorem{continua}[th1]{Theorem}
\newtheorem{isoarc}[th1]{Theorem}
\newtheorem{compinv1}[th1]{Theorem}
\newtheorem{compinv2}[th1]{Theorem}
\newtheorem{la1}{Lemma}
\newtheorem{la2}[la1]{Lemma}
\newtheorem{la3}[la1]{Lemma}
\newtheorem{la4}[la1]{Lemma}

\newtheorem{ahlf}[la1]{Lemma}
\newtheorem{hypmetric}[la1]{Lemma}
\newtheorem{classB}[la1]{Lemma}
\newtheorem{bky}[la1]{Lemma}
\newtheorem{co1}{Corollary}
\newtheorem{co2}[co1]{Corollary}
\newtheorem{co3}[co1]{Corollary}
\newtheorem{co4}[co1]{Corollary}
\theoremstyle{definition}
\newtheorem{qu1}{Question}
\newtheorem{qu2}[qu1]{Question}
\newtheorem{qu3}[qu1]{Question}
\newtheorem{qu4}[qu1]{Question}
\newtheorem{qu5}[qu1]{Question}
\newtheorem{qu6}[qu1]{Question}
\newtheorem{qu7}[qu1]{Question}
\newtheorem{qu8}[qu1]{Question}
\newtheorem{qu9}[qu1]{Question}
\newtheorem{qu10}[qu1]{Question}
\newtheorem{qu11}[qu1]{Question}
\newtheorem{qu12}[qu1]{Question}
\newtheorem{qu13}[qu1]{Question}
\newtheorem{qu14}[qu1]{Question}
\newtheorem{qu15}[qu1]{Question}
\newtheorem{qu16}[qu1]{Question}
\newtheorem{qu17}[qu1]{Question}
\newtheorem{qu18}[qu1]{Question}
\newtheorem{qu19}[qu1]{Question}
  
\newcommand{\C}{{\Bbb C}}
\newcommand{\R}{{\Bbb R}}
\newcommand{\N}{{\Bbb N}}
\newcommand{\Z}{{\Bbb Z}}
\newcommand{\Q}{{\Bbb Q}}
\newcommand{\Si}{\operatorname{sing}}
\newcommand{\re}{\operatorname{Re}}
\newcommand{\im}{\operatorname{Im}}
\newcommand{\Int}{\operatorname{int}}
\newcommand{\Ann}{\operatorname{ann}}
\newcommand{\Deg}{\operatorname{deg}}
\newcommand{\Arg}{\operatorname{arg}}
\newcommand{\Meas}{\operatorname{meas}}

\begin{document}
\def\currentvolume{29}
\def\currentissue{2}
\def\currentyear{1993}
\def\currentmonth{October}
\def\copyrightyear{1993}
\def\currentpages{151-188}
\title{Iteration of meromorphic functions}
\author{Walter Bergweiler}
\address{Lehrstuhl II f\"ur Mathematik, RWTH Aachen,
D-52056 Aachen, Germany}
\email{sf010be@@dacth11.bitnet}
\date{February 23, 1993 and, in revised form, May 1, 1993}
\keywords{Iteration, meromorphic function, entire function,
set of normality, Fatou~set, Julia set, periodic point,
wandering domain, Baker domain, Newton's method}
\subjclass{Primary 30D05, 58F08; Secondary 30D30, 65H05}
\maketitle
 
\section*{Contents}
\par
\noindent
1.\quad Introduction\newline

\noindent
2.\quad Fatou and Julia Sets\newline
\phantom{2.\quad}2.1.\ The definition of Fatou and Julia 
sets\newline
\phantom{2.\quad}2.2.\ Elementary properties of Fatou and 
Julia sets\newline

\noindent
3.\quad Periodic Points\newline
\phantom{2.\quad}3.1.\ Definitions\newline
\phantom{2.\quad}3.2.\ Existence of periodic points\newline
\phantom{2.\quad}3.3.\ The Julia set is perfect\newline
\phantom{2.\quad}3.4.\ Julia's approach\newline

\noindent 4.\quad The Components of the Fatou set\newline
\phantom{2.\quad}4.1.\ The types of domains of 
normality\newline
\phantom{2.\quad}4.2.\ The classification of periodic 
components\newline
\phantom{2.\quad}4.3.\ The role of the singularities of 
the inverse
function\newline
\phantom{2.\quad}4.4.\ The connectivity of the components 
of the Fatou set\newline
\phantom{2.\quad}4.5.\ Wandering domains\newline
\phantom{2.\quad}4.6.\ Classes of functions without 
wandering domains\newline
\phantom{2.\quad}4.7.\ Baker domains\newline
\phantom{2.\quad}4.8.\ Classes of functions without Baker 
domains\newline
\phantom{2.\quad}4.9.\ Completely invariant domains\newline

\noindent
5.\quad Properties of the Julia Set\newline
\phantom{2.\quad}5.1.\ Cantor sets and real Julia 
sets\newline
\phantom{2.\quad}5.2.\ Points that tend to infinity\newline
\phantom{2.\quad}5.3.\ Cantor bouquets\newline

\noindent
6.\quad Newton's Method\newline
\phantom{2.\quad}6.1.\ The unrelaxed Newton method\newline
\phantom{2.\quad}6.2.\ The relaxed Newton method\newline

\noindent
7.\quad Miscellaneous\ topics\newline

\noindent
References

\section{Introduction} \label{intro}
Mathematical models for phenomena in the natural sciences 
often lead
to iteration. An often-quoted example (compare \cite{May76})
comes from population biology. Assuming that the size of a
generation of a population depends solely on the size
of the previous generation and may thus be expressed as a 
function
of it, questions concerning the further development of the 
population
reduce to iteration of this function.
More often, a phenomenon from physics or other sciences
is described by a differential equation.
In certain cases,
for example, if there is a periodic solution,
this differential equation may be studied by
looking at a Poincar\'e return map (see, e.g., \cite[\S 
1.4]{Rue89}),
and again we are led to iteration.
If we solve the differential equation numerically, we are 
also
likely to use a method based on iteration. In fact, many 
algorithms
of numerical analysis (not only those for solving 
differential
equations) involve iteration. One such algorithm, Newton's 
method
of finding zeros,
will be discussed in some detail in \S \ref{newton}.
Apart from that section, however, we will mainly study 
iteration
theory in its own right without having specific applications
in mind.
On the other hand, it is hoped that the questions 
considered here
may also serve as models for other situations so that
their study will enhance our knowledge
of dynamical systems in general.
 
There are two basic problems in iteration theory. The first
(and classical) one is
to study the iterative behavior of an individual function; 
the second
one is to study how the behavior changes if the function is
perturbed, the simplest (but already sufficiently 
complicated)
case being a family of functions that depends on one 
parameter.
Although the second aspect has received much attention
in recent years, we shall consider here only the first one,
except for a few short remarks in \S \ref{misc}.
On the other hand, a good understanding of the dynamics of
an individual function is of course necessary for the study
of problems involving perturbation of functions.
We shall restrict ourselves to functions of one complex 
variable
that are meromorphic in the complex plane. This includes 
rational
and entire functions as special cases.
 
Although some work on iteration was already done in the 
last century,
it is fair to say that the iteration theory of rational 
functions
originated with the work of Fatou \cite{Fat19} and Julia 
\cite{Jul18},
who published long memoirs on the subject between 1918 and 
1920.
At least Fatou's motivation was partly to study functional 
equations,
yet another reason to consider iteration theory.
At the same time, the iteration of rational functions
was also investigated by Ritt \cite{Rit20}.
Some years later in 1926 Fatou \cite{Fat26} extended some of
the results to the case of transcendental entire functions.
He did not, however, consider transcendental meromorphic 
functions, because
\cite[p.\ 337]{Fat26}
here, in general, the iterates have infinitely many 
essential singularities.
Julia did not consider the iteration of transcendental 
functions at all.
(As already pointed out by Fatou \cite[p.\ 358]{Fat26}, 
there occurs a
serious difficulty when one tries to generalize Julia's 
approach to
the transcendental case; see the discussion in \S 
\ref{julia}.)
 
In the past decade there was a renewed interest in the 
iteration theory
of analytic functions, partially due to the beautiful 
computer graphics
related to it (see, for example, the book by Peitgen and 
Richter
\cite{Pei86}),
partially due to new and powerful mathematical methods
introduced into it (notably those introduced by
Douady and Hubbard \cite{Dou82} and by
Sullivan \cite{Sul82}).
Most of the work has centered around the iteration
theory of rational functions,
but there is also a considerable number of
papers devoted to transcendental entire functions; and in 
recent years
work on the iteration of transcendental meromorphic 
functions has also
begun.
 
There exist a number of introductions to and surveys of
the iteration theory of rational functions.
We mention
\cite{Bea91,Bla84,Car90,Dev89a,Dou84,Ere90,Kee89a,Lyu86,
Mil90,Ste92}
among the more recent ones but also some older ones
\cite{Bro67,Cre25,Mon27}.
There are comparatively few expositions of the iteration 
theory
of transcendental functions. We refer to \cite{Bak87}
for the iteration of entire functions and to \cite{Ere90}, 
which has
a chapter on this topic.
 
This paper attempts to describe some of the results obtained
in the
iteration theory of transcendental meromorphic functions,
not
excluding the case of entire functions. The reader is not
expected to be familiar with the iteration theory
of rational functions.
On the other hand,
some aspects where the transcendental case is
analogous to the rational case are treated rather briefly 
here.
For example, we introduce the different types of components
of the Fatou set that occur in the iteration of rational
functions but omit a detailed description of these types.
Instead,  we concentrate on the types of components
that are special to transcendental functions
(Baker domains and wandering domains).
 
This article is mainly an exposition of known results, but
it also contains
some new results. For example, Theorems \ref{th2} and 
\ref{bakdom} have
been known before only for entire functions or special 
classes of
meromorphic functions.
Other results like Theorem \ref{th7} or Corollaries 
\ref{co1}
or \ref{co2} are certainly known to those who work in 
the field,
but they do not seem to have been stated explicitly 
before.
 
As already mentioned, there are beautiful computer graphics
related
to the theory, and there are many places (besides 
\cite{Pei86})
where such pictures can be found for rational functions.
Although Julia sets (and bifurcation diagrams)
of transcendental functions
can compete in their beauty and complexity  very well with
those of rational functions,
this article is not illustrated with  such pictures.
The interested reader is referred to
[53--56, 96].

\section{Fatou and Julia sets} \label{fatou}
\subsection{The definition of Fatou and Julia sets}  
\label{deffatou}
Let $f:\C\to\widehat{\C}$ be a meromorphic function, 
where $\C$ is the
complex plane and $\widehat{\C}=\C\cup\{\infty\}$.
Throughout this paper, we shall always assume that $f$ is 
neither constant
nor a linear transformation.
Denote by $f^n$
the $n$th iterate of $f$, that is, $f^{\,0}(z)=z$ and 
$f^n(z)=f(f^{n-1}(z))$
for $n\geq 1$. Then $f^n(z)$ is defined for all $z\in\C$ 
except for
a countable set which consists of the poles of $f,f^2,
\dots,f^{n-1}$.
If $f$ is rational, then $f$ has a meromorphic extension 
to $\widehat{\C}$;
and, denoting the extension again by $f$,
we see that $f^n$ is defined and meromorphic in 
$\widehat{\C}$.
But if $f$ is transcendental---and this is the case we are 
mainly
interested in---there is, of course,
no (reasonable) way to define $f(\infty)$.
 
The basic objects studied in iteration theory are the
Fatou set $F=F(f)$ and the Julia set $J=J(f)$
of a meromorphic function $f$.
Roughly speaking, the Fatou set is the set where the 
iterative
behavior is relatively tame in the sense that points
close to each other behave similarly,
while the Julia set is the set where chaotic phenomena take 
place.
The formal definitions are
\[
F=\{z\in\widehat{\C}: \{f^n: n\in \N\}
\text{\ is defined and normal in some neighborhood of\ } z\}
\]
and
\[
J=\widehat{\C}\backslash F.
\]
As already mentioned, the requirement that $f^n$ be defined 
is always
satisfied if $f$ is rational, and hence it can be (and of 
course
always is) omitted from the definition.
An analogous remark applies to transcendental
entire functions, where $f^n$ is defined
for all $z\in\C$. In this case, we always have $\infty\in 
J$.
 
A similar case is given by meromorphic functions with
exactly 
one pole
if this pole is  an omitted value.
(A complex number $z_0$ is called an {\em omitted value} of 
the
meromorphic function $f$,  if $f(z)\neq z_0$ for all $z\in 
\C$.)
In this case, if the pole of $f$ is denoted by $z_0$, we 
have
$\{z_0,\infty\}\subset J$, and $f^n(z)$ is defined for all
$z\in \widehat{\C}\backslash \{z_0,\infty\}$.
 
It is not difficult to show that $f$ has the form
\[
f(z)=z_0+\frac{e^{g(z)}}{(z-z_0)^m}
\]
for some positive integer $m$ and some entire function $g$
in 
this case.
It is no loss of generality to assume that $z_0=0$
so that
\begin{equation}
f(z)=\frac{e^{g(z)}}{z^m},  \label{P}
\end{equation}
because otherwise we may consider $\phi^{-1}(f(\phi(z)))$
instead of $f$,
where $\phi(z)=z+z_0$. More generally, instead of maps of 
the form
(\ref{P}), we may consider
analytic self-maps of $\C\backslash \{0\}$ here, without 
requiring that
$0$ be a pole of the map; compare
\cite{Bak87a,Bha69,Kee88,Kee89,Kot87,Kot90,Mak88,Rad53}.
We shall restrict
ourselves, however, to the case where $f$ is meromorphic in 
$\C$,
our main interest being in the case where $f$ is entire or 
has several poles anyway.
 
In the remaining case, where $f$ has either at least two 
poles
or only one pole which is not an omitted value,
there are infinitely many points that are mapped onto a pole
of $f$ by some iterate of $f$.
For $z_0\in\widehat{\C}$, we define the {\em backward orbit}
$O^-(z_0)$ of $z_0$ by
\[
O^-(z_0)=\bigcup_{n\geq 0}f^{-n}(z_0),
\]
where $f^{-n}(z_0)=\{z:f^n(z)=z_0\}$.
Then the above statement is equivalent to saying that
$O^-(\infty)$ is an infinite set.
In fact, already $f^{-3}(\infty)$ is infinite, as follows
easily
from Picard's theorem. The largest open set where all 
iterates
are defined is given by $\widehat{\C}\backslash 
\overline{O^-(\infty)}$.
Since $f(\widehat{\C}\backslash \overline{O^-(\infty)})
\subset \widehat{\C}\backslash \overline{O^-(\infty)}$
and since $O^-(\infty)$ has more than two elements, 
$\{f^n\}$ is normal in
$\widehat{\C}\backslash \overline{O^-(\infty)}$ by Montel's 
theorem.
Hence
\[
F=\widehat{\C}\backslash \overline{O^-(\infty)}
\quad\text{and}\quad
J=\overline{O^-(\infty)};
\]
compare \cite{Bak91} and \cite{Rad53}.
We see that in this case the requirement that $\{f^n\}$ be
normal
can be omitted from the definition.
 
From this point of view, the iteration theory of entire 
functions
and of meromorphic functions with one pole which is an 
omitted value
is quite different from that of general meromorphic 
functions, which
have at least two poles or only one pole which is not 
omitted.
In the first two cases, it is clear where the iterates are
 defined,
and we ask where they form a normal family. In the third 
case,
we just ask where they are defined, and this then implies 
that they
form a normal family there. Therefore, it is not surprising 
that
there are major differences between these cases. On the 
other hand,
there are also many analogies, but even here the proofs are 
often
quite different.
 
According to the above remarks we shall divide the class of
transcendental meromorphic
functions for further reference into three subclasses:
\begin{itemize}
\item $E=\{f: f \makebox{\ is transcendental entire}\};$
\item $P=\{f: f$ is transcendental meromorphic,
 has exactly one pole,\newline
\qquad\qquad\quad and this pole
is an omitted value$\}$;
\item $M=\{f: f$ is transcendental meromorphic
 and has either at least\newline
\qquad \qquad\quad two poles or
exactly  one pole which is not an omitted value$\}$.
\end{itemize}
Here $E$ and $M$ are thought of as mnemonics
for {\em e}ntire and (general) {\em m}eromorphic functions, 
while
$P$ stands for one {\em p}ole (or {\em p}unctured plane).
As already mentioned, we may (and often will) assume that 
functions
in $P$ have the form (\ref{P}).
 
\subsection{Elementary properties of Fatou and Julia sets} 
\label{elem}
By definition, $F$ is open and $J$ is closed. The
properties of $F$ and $J$ contained in the following
three lemmas are also easily verified,
the proofs for transcendental functions being analogous to 
those
for the rational case.
\begin{la1}  \label{la1}
If $f$ is rational, $f\in P$, or $f\in E$,
then $F(f)=F(f^n)$ and $J(f)=J(f^n)$ for all $n\geq 2$.
\end{la1}
Here we have to exclude $f\in M$, because then $f^n$ is not 
meromorphic
in $\C$ so that $F(f^n)$ and $J(f^n)$ are not defined.
(There is, of course, a natural way to define $F(f^n)$ for 
$f\in M$
and $n\geq 2$, or, more generally, to define $F(f)$ for 
functions $f$
meromorphic in $\C$ except for countably many points.
Then the conclusion of Lemma \ref{la1} holds for such 
functions.)
\begin{la2} \label{la2}
$F$ and $J$ are completely invariant.
\end{la2}
Here, by definition, a set $S$ is
called completely invariant if $z\in S$ implies that
$f(z)\in S$, unless $f(z)$ is undefined, and
that $w\in S$ for all $w$ satisfying $f(w)=z$.
\begin{la3} \label{la3}
Either $J=\widehat{\C}$ or $J$ has empty interior.
\end{la3}
We note that the case $J=\widehat{\C}$ is actually possible.
Examples of rational functions with this property are
given by the rational functions that come from the 
multiplication
theorems of elliptic functions.
Usually, Latt\`es (\cite{Lat18}, see also \cite{Jul18a})
is credited with having introduced them
into the subject, but it should perhaps be mentioned that
already B\"ottcher \cite[p.\ 63]{Boe98} was aware of these 
examples.
For other examples of rational functions satisfying
 $J=\widehat{\C}$
we refer to \cite[\S 9.4, \S 11.9]{Bea91}.
The first example of an entire function with this property
was given by Baker \cite{Bak70}, who proved that
$J(\lambda ze^z)=\widehat{\C}$ for a suitable value of 
$\lambda$.
Later, Misiurewicz \cite{Mis81} proved that 
$J(e^z)=\widehat{\C}$,
confirming a conjecture of Fatou \cite[p.\ 370]{Fat26}.
To obtain an example in $P$, we note that modifications
of Baker's argument show that 
$J(\lambda e^z/z)=\widehat{\C}$ for a suitable
value of $\lambda$. Finally, as an example in $M$, we 
mention
that $J(\lambda\tan z)=\widehat{\C}$ for suitable values of 
$\lambda$.
We comment on these examples in \S \ref{nobakdom}.
We note that while $J=\widehat{\C}$ is possible, we always
have $F\neq \widehat{\C}$. In fact, as we shall see in
\S \ref{perfect}, $J$ is a perfect (and hence uncountable) 
set.
 
We say that $z_0$ is {\em exceptional} if $O^-(z_0)$ is 
finite.
It is not difficult to see that
meromorphic functions have at most two exceptional values.
(For transcendental functions this is an immediate 
consequence
of Picard's theorem.) For rational functions the exceptional
values
are in the Fatou set.
This is not necessarily the case for transcendental 
meromorphic functions.
If $f\in P$ with pole at $z_0$, then
$z_0$ and $\infty$ are exceptional, but
there are no further exceptional values.
Similarly,
if $f\in E$, then $\infty$ is always exceptional so
that there is at most one finite exceptional value.
If a function $f\in E$ has a finite exceptional value $z_0$,
then it has the form $f(z)=z_0+(z-z_0)^me^{g(z)}$ for
some nonnegative integer $m$ and an entire function $g$ so
that $f$ can also be considered as a self-map of the 
punctured
plane.
 
A simple consequence of Montel's theorem is the following 
result.
\begin{la4} \label{la4}
If $z_0\in J$ is not exceptional, then 
$J=\overline{O^-(z_0)}$.
\end{la4}
Similar to the backward orbit $O^-(z_0)$, we
define the {\em forward orbit}
$O^+(z_0)$ of $z_0\in\widehat{\C}$ by
$O^+(z_0)=\bigcup_{n\geq 0}f^n(z_0)$. Of course,
here the union is taken only over those $n\geq 0$ for which
$f^n(z_0)$ is defined.
The {\em orbit}
$O(z_0)$ of $z_0$ is defined by 
$O(z_0)=O^+(z_0)\cup O^-(z_0)$.
For a subset $S$ of $\widehat{\C}$,
we put $O^{\pm}(S)=\bigcup_{z\in  S}O^\pm(z)$
and
$O(S)=\bigcup_{z\in  S}O(z)$.
With this terminology, Lemma \ref{la2}
may be written in the form $O(F)\subset F$ and 
$O(J)\subset J$.
 
Another simple consequence of Montel's theorem is that if 
$U$ is an open
set that contains a point of $J$, then 
$\widehat{\C}\backslash O^+(U)$
contains at most two points, and these points are 
exceptional.
For rational $f$, this implies that
$O^+(J\cap U)=J$ and even that $f^n(J\cap U)=J$
for all sufficiently large $n$; see 
\cite[Theorem 4.2.5]{Bea91}.
For transcendental $f$ we find that 
$J\backslash O^+(J\cap U)$
contains at most two points and that these points are 
exceptional.
However, we cannot
deduce from this that
$J\backslash f^n(J\cap U)$ contains  only two points for 
sufficiently
large $n$.
In general,
$J\backslash f^n(J\cap U)$ will contain neighborhoods
of the exceptional points (if they exist).
But if $f\in M$ does not have exceptional points or if they 
are in $F$,
then $f^n(J\cap U)=J$
for all sufficiently large~$n$.
 
\section{Periodic points}   \label{period}
\subsection{Definitions} \label{defperiod}
An important role in iteration theory is played by the 
periodic
points. By definition, $z_0$ is called a {\em periodic 
point}
of $f$ if $f^n(z_0)=z_0$ for some $n\geq 1$. In this case, 
$n$
is called a {\em period} of $z_0$, and the smallest $n$ with
this
property is called the {\em minimal period} of $z_0$.
For a periodic
point $z_0$ of minimal period $n$,
$(f^n)'(z_0)$ is called the {\em multiplier} of $z_0$.
(If $z_0=\infty$, which can happen only for rational
function $f$, of course, this has to be modified.
In this case, the multiplier is defined to be
$(g^n)'(0)$ where $g(z)=1/f(1/z)$.)
A periodic point is called {\em attracting}, 
{\em indifferent}, or
{\em repelling\/} accordingly as
the modulus of its multiplier
is less than, equal to, or greater than $1$.
Periodic points of multiplier $0$ are called
{\em superattracting}.
(Some writers reserve the term {\em attracting\/} for
the case
$0<|(f^n)'(z_0)|<1$, but we consider superattracting as a
 special case
of attracting.)
The multiplier of an indifferent
periodic point is of the form
$e^{2\pi i \alpha}$ where $0\leq \alpha <1$.
We say that $z_0$ is {\em rationally indifferent} if
$\alpha$ is rational and {\em irrationally indifferent}
otherwise. Also, a point $z_0$ is called {\em preperiodic}
if $f^n(z_0)$ is periodic for some $n\geq 1$.
Finally, a periodic point of period 1 is called
a {\em fixed point}.
 
It is easy to see that attracting periodic points are in
 $F$,
while repelling and rationally indifferent periodic points
are in $J$.
For irrationally indifferent periodic points the question
 whether
they are in $F$ or $J$ is difficult to decide.
Both possibilities do occur. We refer the reader to the
 classical
papers of Cremer \cite{Cre28,Cre38} and Siegel \cite{Sie42},
as well as more recent work of Yoccoz \cite{Yoc87,Yoc88}.
An exposition of these  and other results, together with
 further references,
can be found in \cite{Per92}.
 
The behavior of the iterates in the neighborhood of a fixed 
point
(or, more generally, a periodic point)
is intimately connected with the solution of certain 
functional
equations. For most results in this direction,
 it is required only that the function
under consideration
is defined in a neighborhood of the fixed point, and
it is usually irrelevant whether it extends to a rational or
 transcendental
meromorphic function.
Therefore, we omit this topic here but refer to the papers 
and books
on iteration of rational functions cited in the 
introduction.
 
\subsection{Existence of periodic points} \label{existence}
It is clear that a rational function has periodic points of
(not necessarily minimal) period
$n$ for all $n\geq 1$. Transcendental entire functions, 
however,
need not have fixed points, that is, periodic points of 
period $1$.
A simple example is given by $f(z)=e^z+z$.
On the other hand, already Fatou
\cite[p.\ 345]{Fat26} proved that an entire transcendental
function $f$ has at least one periodic point of period 2.
The idea is to consider the function
\[
h(z)=\frac{f(f(z))-z}{f(z)-z}.
\]
If $f$ does not have periodic points of period $2$ (and 
hence
does not have fixed points), then $h$ is an entire function 
which does
not take the values $0$ and $1$. By Picard's theorem $h$ is
constant. Once this is known, it is not difficult to obtain 
a
contradiction.
Fatou's result was generalized by
Rosenbloom \cite{Ros48}, who proved
the following theorem.
\begin{th1} \label{th1}
An  entire transcendental
function has infinitely many
periodic points of period $n$ for all $n\geq 2$.
\end{th1}
The idea of  the proof is similar.  Instead of Picard's 
theorem,
however, Rosenbloom used something stronger,
namely, Nevanlinna's theory  on the distribution
of values, which may be considered as a quantitative version
of Picard's theorem. Since this is the only place in this 
paper where we
use Nevanlinna theory, we do not give an introduction to it
but refer to \cite{Hay64,Jan85,Nev70} for notation and basic
results.
 
To prove Theorem \ref{th1}, we suppose that $f^n$ and hence
$f$ have only finitely many fixed points and
consider the auxiliary function
$h$ defined by
\[
h(z)=\frac{f^n(z)-z}{f^{n-1}(z)-z}.
\]
Then
\begin{equation}
\begin{split}
N(r,h)&=O(T(r,f^{n-1})),\\
N\left(r,\frac{1}{h}\right)&=O(\log r),
\end{split}\nonumber
\end{equation}
and
\[
N\left(r,\frac{1}{h-1}\right)=O(T(r,f^{n-1}))
\]
as $r\to\infty$.
Also, it is not difficult to prove that 
$T(r,f^{n-1})=o(T(r,f^n))$
as $r\to\infty$ outside some exceptional set
of finite measure; see \cite[p.\ 147]{Jan85}.
(In fact, this last result even holds without exceptional 
set;
see \cite[Theorems 1 and 2]{Clu70} or 
\cite[Lemma 2.6]{Hay64}.)
We deduce that $T(r,h)\sim T(r,f^n)$ and hence
that
\[
N(r,h)+N\left(r,\frac{1}{h}\right)+N\left(r,\frac{1}{h-1}
\right)=o(T(r,h))
\]
as $r\to\infty$
outside the exceptional set.
This contradicts Nevanlinna's second fundamental theorem.
Thus the proof of Theorem \ref{th1} is  complete.
 
One may ask whether there is some quantitative version of 
Theorem
\ref{th1} in the sense that there is a lower bound for 
$N(r,1/(f^n(z)-z))$
in terms of $T(r,f^n)$ if $n\geq 2$ and $f\in E$.
Denote by $\delta(a,h)$ the deficiency of a meromorphic $h$
with respect to the value $a\in\widehat{\C}$.
\begin{qu1}
Do we have
$\delta(0,f^n(z)-z)=0$ (or at least
$\delta(0,f^n(z)-z)<1$) if $f\in E$ and $n\geq 2$?
\end{qu1}
Some (much weaker) results of this
type can be found in \cite{Bak59,Bak59-60,Ber90,Yan91}.
 
If $f\in P$, then the conclusion of Theorem \ref{th1}
is true even if $n=1$. To see this,
suppose $f\in P$ and define $F(z)=f(z)/z$. Then
$F$ has no zeros and
at most two poles.
Hence,
by Picard's theorem $F$ takes the value $1$ infinitely 
often,
that is, $f$ has infinitely many fixed points.
If $f\in M$, then $f$ need not have fixed points;
consider $f(z)=z+1/g(z)$ where $g$ is entire transcendental.
But we shall see below that $f$ has infinitely many periodic
 points
of period $n$ if $n\geq 2$.
 
Baker [4--6,8]
seems to have been the first
who addressed the question in which cases a
rational or entire function may fail to
have periodic points of minimal period $n$ for some $n$.
He proved \cite{Bak64} that
if a rational function $f$ of degree $d\geq 2$ has no 
periodic point
of minimal period $n$, then
$n=2$ and $d\in \{2,3,4\}$ or $n=3$ and $d=2$.
Moreover, if $f$ is a polynomial, then only the case $n=d=2$
can
occur.
Earlier he had
proved \cite{Bak60} that if $f$ is an entire function,
then there exists at most one integer $n\geq 1$ (depending 
on $f$)
such
that $f$ does not have periodic points of minimal 
period~$n$.
 
The latter result can be strengthened as follows.
 
\begin{th1a} \label{th1a}
If $f$ is a transcendental meromorphic function and 
$n\geq 2$,
then $f$ has infinitely many periodic points of minimal 
period $n$.
\end{th1a}
As already mentioned, this also holds for $n=1$ if $f\in P$,
but not
in general if $f\in E$ or $f\in M$.
 
Theorem \ref{th1a} was proved in \cite{Ber91}
if $f\in E$ and in \cite[Chapter 5.2]{Bha69} if $f\in P$.
A proof for $f\in M$, using the ideas of \cite{Ber92},
is as follows.
 
Suppose first that $f^{n-1}$ has at least three poles 
$p_1,p_2,p_3$.
Define $g=f^{n-1}$ and denote by $m_j$ the order of the pole
$p_j$.
There exist functions $h_j$, defined and analytic in a 
neighborhood
of $0$, such that $g(p_j+h_j(z))=z^{-m_j}$.
Define $k_1(z)=p_1+h_1(z^{m_2m_3})$,
$k_2(z)=p_2+h_2(z^{m_1m_3})$,
and $k_3(z)=p_3+h_3(z^{m_1m_2})$ so that 
$g(k_j(z))=z^{-m_1m_2m_3}$.
Suppose now that $f$ does not have periodic points of period
$n$ in
neighborhoods of $p_1,p_2,p_3$.
Then $F(z)=f(z^{-m_1m_2m_3})=f^n(k_j(z))\neq k_j(z)$ in a
neighborhood of $0$. Hence
\[
\frac{(F(z)-k_1(z))(k_3(z)-k_2(z))}
{(F(z)-k_2(z))(k_3(z)-k_1(z))}\neq 0,1,\infty
\]
in a neighborhood of $0$, contradicting Picard's theorem.
Hence the periodic points of period $n$ must accumulate at 
$p_1$,
$p_2$, or $p_3$. Since periodic points of period $k$ where 
$1\leq k \leq n-1$
cannot accumulate at poles of $f^{n-1}$, we deduce that $f$ 
has infinitely
many periodic points of minimal period $n$.
 
The case that $f^{n-1}$ has exactly two poles is similar.
Here we choose $k_1(z)=p_1+h_1(z^{m_2})$,
$k_2(z)=p_2+h_2(z^{m_1})$, $F(z)=f(z^{-m_1m_2})$ and
observe that
\[
\frac{F(z)-k_1(z)}{F(z)-k_2(z)}\neq 0,1,\infty
\]
in a neighborhood of $0$.
 
Finally, we
consider the case that $f^{n-1}$ has only one pole but is in
$M$.
Then necessarily $n=2$ or $n=3$.
We leave it to the reader to check that the method
of \cite{Ber91} can be extended to this case.
 
It is perhaps worth mentioning
(and, at least for me, surprising)
that the proof given above for functions
with at least two poles is shorter and more elementary than 
the proofs
for the classes $E$ and $P$ so that
these questions are much
simpler for functions with poles than for entire functions.
We will give a generalization of
Theorem \ref{th1a}
in \S \ref{julia}. Its proof, however, will be less
elementary (but still short).
 
\subsection{The Julia set is perfect}  \label{perfect}
The results concerning the existence of periodic points may
be used to prove that $J\neq \emptyset$.
More generally, we have the following result.
\begin{perf}
Let $f$ be a meromorphic function. Then $J(f)$ is perfect.
\end{perf}
Recall that a set is called perfect if it is closed,
nonempty, and does not contain isolated points.
 
We first prove that $J$ is not empty and in fact an infinite
set.
There are essentially two ways to
do this if $f$ is rational. One
method is to assume that $f^{n_j}\to\phi$ uniformly in 
$\widehat{\C}$.
Then $\phi$ must also be rational and
$\Deg(\phi)=\lim_{j\to\infty} \Deg(f^{n_j})$. But
$\Deg(f^{n_j})=(\Deg(f))^{n_j}\to\infty$ as $j\to\infty$,
provided $\Deg(f)\geq 2$, a contradiction.
The other method (which is the one used by Fatou and Julia)
is to prove that $f$ has a
fixed point which is  repelling or has multiplier $1$.
Once  a point $z_0\in J$ is found, it is not difficult to 
see that
$O^-(z_0)$ is infinite.  Hence $J$ is infinite because
$O^-(z_0)\subset J$.
 
Once this is known, we can prove that
$J$ is perfect as follows. Suppose that $w_0\in J$, and let
$N$ be a neighborhood of $w_0$. We can find $w_1,w_2,w_3\in
J\backslash O^+(w_0)$. Because $\{f^n|_N\}$ is not normal,
\!$w_j\!\in\! O^+(N)$ for some 
\!$j\!\in\!\{1\!,\!2\!,\!3\}$\!. Hence
$O^-(w_j)\cap N\backslash\{w_0\}$ is not empty.
In particular, $J\cap N\backslash \{w_0\}$ is not empty.
Hence $w_0$ is not isolated, that is, $J$ is perfect.
 
Both methods to prove that $J$ is infinite
do not generalize to the case where $f$ is transcendental
(see, however, the discussion of the second method in the 
next section).
On the other hand, it is clear  from the discussion in
\S \ref{deffatou} that $J(f)$ is an infinite and, in fact,
a perfect set if $f\in M$.
 
We sketch how one can prove that $J(f)$ is infinite in the 
case that
$f\in E$ or $f\in P$.
First we note that we may assume $f$ has infinitely many
fixed points. As already mentioned above, this is always the
 case
if $f\in P$. And if $f\in E$, then
we may consider $f^{\,2}$ instead of
$f$  because of Lemma \ref{la1}, and $f^{\,2}$ always has 
infinitely
many fixed points by Theorem \ref{th1}.
If infinitely many of the fixed points of $f$ are in $J$, 
then
we are done. Hence we may assume that there exist two fixed
points $p$ and $q$ of $f$ that are contained in $F$.
If $p$ and $q$ are in different components of $F$, then any 
path
connecting them must meet $J$, and we are also done.
Thus we may assume that there exists a component $U$
of $F$ which contains $p$ and $q$. Clearly, the limit 
functions
of $\{f^n|_U\}$ cannot be constant. We deduce that if
$f^{n_j}(z)\to\phi(z)$ for $z\in U$ as $j\to\infty$,
then $f^{n_{j+1}-n_j}(z)\to z$ for $z\in U$.
This implies that $f|_U$ is an automorphism of $U$.
Hence $f^{-1}(U)$ contains components of $F$ different from 
$U$.
Again, any path connecting these components meets $J$.
This completes the proof that $J$ is an infinite set.
(With a little more effort one can show that
a component of $F$ cannot contain two fixed points; see
\cite[Lemma 6.9.3]{Bea91}, and compare also
Theorem \ref{classification} in \S \ref{periodiccom}.)
The proof that $J$ is perfect can now be carried out
as in the rational case.

\subsection{Julia's approach} \label{julia}
So far our development of the theory has followed Fatou's 
ideas.
Julia based his theory on the closure of the set of 
repelling periodic
points. One of the basic results of the theory is that these
two sets
are actually equal.
\begin{rep} \label{rep}
Let $f$ be a meromorphic function. Then $J(f)$ is the 
closure of the set
of repelling periodic points of $f$.
\end{rep}
For rational $f$, this result was obtained by both Fatou
\cite[\S 30, p.\ 69]{Fat19} and Julia
\cite[p.\ 99, p.\ 118]{Jul18}.
Their proofs, however, were different. (A good exposition of
both
proofs can be found in \cite[\S 11]{Mil90}.)
Fatou proved first that any point in $J$ is the limit point 
of
periodic points and then that there are only finitely many
nonrepelling periodic points, which together implies the 
result.
The first part does carry over to transcendental functions
(\cite[p.\ 354]{Fat26}, see also \cite{Chu85}),
but the second part clearly does not, as can  be seen by 
simple  examples
like $f(z)=e^z+z+1$.
For Julia's method it is essential that the set of repelling
periodic points is not empty. In fact, it suffices that
the set of repelling and rationally indifferent
periodic points is not empty, and
a rational function always has at least one
fixed point which is repelling or has multiplier $1$
(see \cite[\S 2, p.\ 168]{Fat19}
and \cite[p.\ 85, p.\ 243]{Jul18}). Julia's method does 
carry over
to transcendental functions with a repelling
or rationally indifferent
fixed point which is not exceptional (this is done in
\cite[pp.\ 229--230; 59, pp.\ 69--70]{Dev91},
but in general transcendental functions need not have
such a fixed point; in fact, they need not have
fixed points at all.
 
Baker \cite{Bak68} proved that Theorem \ref{rep}
holds for entire transcendental functions as well.
His proof was based on a deep theorem of Ahlfors
\cite{Ahl35,Hay64,Nev70} from his theory of covering 
surfaces.
Theorem \ref{rep} was extended to class $P$ in 
\cite[Theorem 5.2]{Bha69}
and to class $M$ in \cite[Theorem 1]{Bak91},
the proofs being based again on Ahlfors's theorem.
 
We sketch the argument for class $M$ and begin with the 
statement
of a version of Ahlfors's theorem.
(A different version is used in \cite{Bak68} and 
\cite{Bha69}.)
Unfortunately, lack of space prevents us from discussing the
proof
of this important result.
\begin{ahlf}\label{ahlfors}
Let $f$ be a transcendental meromorphic function, and let
$D_1,D_2, \dots,D_5$ be five simply connected domains in 
$\C$ with
disjoint closures.
Then there exists $j\in\{1,2,\dots,5\}$ and, for any
$R>0$, a simply connected
domain $G\subset \{z\in \C:|z|>R\}$ such
that $f$ is a conformal map of $G$ onto $D_j$.
If $f$ has only finitely many poles, then \rom{``}five
\rom{''}
may be replaced by \rom{``}three\rom{''}.
\end{ahlf}
Following \cite{Bak91}, we deduce the following result.
\begin{bky} \label{bky}
Suppose that  $f\in M$ and that
$z_1,z_2,\dots,z_5\in O^-(\infty)\backslash\{\infty\}$
are distinct.
Define $n_j$ by $f^{n_j}(z_j)=\infty$. Then
there exists $j\in\{1,2,\dots,5\}$ such that $z_j$ is
a limit point of repelling periodic points of minimal period
$n_j+1$.
If $f$ has only finitely many poles, then \rom{``}five
\rom{''}
may be replaced by \rom{``}three\rom{''}.
\end{bky}
To deduce Lemma \ref{bky} from Lemma \ref{ahlfors},
we choose the $D_j$  as discs around $z_j$ where the
radii are chosen so small that the $D_j$ do not contain 
critical
points of $f$ and that their closures are pairwise disjoint.
There exists $R>0$ such that
$f^{n_j}(D_j)\supset \{z:|z|>R\}\cup\{\infty\}$.
We choose $j$ and $G$ according to Lemma \ref{ahlfors}.
Then we can find $H\subset D_j$ such that  $f^{n_j}:H\to G$,
and hence $f^{n_j+1}:H\to D_j$ is a conformal mapping. 
Moreover,
$\overline{H}\subset D_j$, and this implies that the inverse
function $f^{-n_j-1}$ of
$f^{n_j+1}:H\to D_j$ has an attracting fixed point in $D_j$.
Clearly, this attracting fixed point of $f^{-n_j-1}$ is a
repelling periodic point of $f$
of period $n_j+1$. Because the $D_j$ can be chosen
arbitrarily small, the repelling periodic
points of period $n_j+1$ accumulate at $z_j$.
Because $z_j$ is a pole of $f^{n_j}$, periodic
points of period less than $n_j+1$ cannot accumulate
at $z_j$; hence, $z_j$ is a limit point
of repelling periodic points of minimal period $n_j+1$.
This completes the proof of Lemma \ref{bky}.
 
Lemma \ref{bky} yields immediately that the conclusion
of Theorem \ref{rep} holds for $f\in M$ because
$J=\overline{O^-(\infty)}$ is
perfect.
 
Another interesting consequence of Lemma \ref{bky} is that
if $f^{-n+1}(\infty)$ contains more than four elements,
then $f$ has infinitely many repelling periodic points of
minimal period $n$. In particular, this is the case if
$f\in M$ and $n\geq 4$.
We also see that
$f$ has infinitely many repelling periodic points of
minimal period $2$ and $3$ if $f$ has more than two poles.
On the other hand, it was proved in \cite{Ber91} that
if $f$ is entire transcendental, then $f$ has infinitely
many repelling periodic points of minimal period $n$
for all $n\geq 2$. The method used there
can be extended to the case where
$n=2$ or $n=3$ and
$f$ has one or two poles. Hence we obtain the following
generalization of Theorem \ref{th1a}.
\begin{th2} \label{th2}
If $f$ is a transcendental meromorphic function and 
$n\geq 2$,
then $f$ has infinitely many repelling
periodic points of minimal period $n$.
\end{th2}
We remark that in view of Lemma \ref{la1},
Julia's method can be used
to obtain Theorem \ref{rep} from  Theorem \ref{th2}. This 
does
not, however,
constitute a new proof of Theorem \ref{rep},
because the argument in \cite{Ber91} also uses Ahlfors's
theorem.
The fact that all known proofs of the existence of repelling
periodic points are based on this deep result makes Julia's
approach to start with the closure of the set of repelling
periodic points inadequate for transcendental functions,
because it is difficult to see that this set is not empty.
 
It would be of interest to give a more elementary proof of
Theorem \ref{rep}.
\begin{qu2}
Is there a proof of Theorem \ref{rep}
which does not use Ahlfors's result?
\end{qu2}
As already mentioned, it suffices to prove the existence of 
a
repelling or rationally indifferent
periodic point which is not exceptional.
 
\section{The components of the Fatou set}   \label{domains}
\subsection{The types of domains of normality}   
\label{types}
Let $U$ be a (maximal) domain of normality of the iterates 
of $f$, that
is, a component of $F$.
(Here and in the following, ``component'' always means 
``connected
component''.)
Then $f^n(U)$ is contained in a component of $F$ which we 
denote
by $U_n$.
A component $U$ is called {\em preperiodic}
if there exist $n>m\geq 0$ such that $U_n=U_m$. In 
particular,
if this is the case for $m=0$ (where $U_0=U$) and some 
$n\geq 1$,
then $U$ is called {\em periodic} with {\em period} $n$,
and $\{U,U_1,\dots,U_{n-1}\}$ is called a (periodic) 
{\em cycle}
of components. Again,
the smallest $n$ with this property is called the
{\em minimal period} of $U$.
In the case $n=1$, that is, if $f(U)\subset U$,
$U$ is called {\em invariant}.
A component of $F$ which is
not preperiodic is called a {\em wandering} component (or
{\em wandering domain}).
 
For rational functions, we have $f(U)=U_1$, but for
transcendental functions it is possible that $f(U)\neq U_1$.
For example, if $f(z)=\lambda e^z$ where $0<\lambda<e^{-1}$,
then
$F$ consists of a single component which contains $0$, but
clearly $0\notin f(F)$. Similarly, if $f(z)=\lambda\tan z$
where $0<\lambda<1$, then
$F$ consists of a single component which contains 
$\pm \lambda i$,
but $\pm\lambda i\notin f(F)$.
 
Values in $U_1\backslash f(U)$ need not be omitted values.
As an example, we consider $f(z)=z\exp((-z^2+3z-2)/6)$.
Then $0$ and $2$ are attracting fixed points, while $1$ is
a repelling fixed point. Let $V$ be the component of $F$
that contains $0$. All large positive real numbers are 
contained
in a component $U$ satisfying $f(U)\subset V$; that is, we 
have
$V=U_1$. It is not difficult to show that $U\neq V$. (For 
instance,
this follows from the fact that $U$ and $V$ are simply 
connected and
symmetric with respect to the real axis.) It follows that
$0\in U_1\backslash f(U)$.
 
On the other hand, it is easy to see that values in
$U_1\backslash f(U)$ are asymptotic values of $f$,
the asymptotic path being contained in $U$. As pointed out 
in
\cite[p.\ 242]{Bak91}, one can deduce from Gross's star 
theorem
that $f(U)$ is a dense open subset of $U_1$.
If $f\in E$, then $U_1\backslash f(U)$
contains at most one point. I.\ N.\ Baker has kindly 
informed me
that this result  was proved by M.\ Herring. It  had also 
been
obtained independently in \cite{Ber93}.
\begin{qu18}
Can $U_1\backslash f(U)$ have more than two points
if  $f\in M$?
\end{qu18}
 
\subsection{The classification of periodic components} 
\label{periodiccom}
The behavior of $f^n$ in periodic components is well
understood.
\begin{classification} \label{classification}
Let $U$ be a periodic component of period $p$.
Then we have one of the following possibilities\/\rom{:}
\begin{itemize}
\item $U$ contains an attracting periodic point $z_0$ of 
period $p$.
Then $f^{np}(z)\to z_0$ for $z\in U$ as $n\to \infty$, and  
$U$
is called the immediate attractive basin of $z_0$.
\item $\partial U$ contains a periodic point $z_0$ of period
$p$ and
$f^{np}(z)\to z_0$ for $z\in U$ as $n\to \infty$.
Then $(f^p)'(z_0)=1$ if $z_0\in\C$. \rom{(}For $z_0=\infty$
we have
$(g^p)'(0)=1$ where $g(z)=1/f(1/z)$.\rom{)}
In this case, $U$ is called a Leau domain.
\item There exists an analytic homeomorphism $\phi:U\to D$  
where
$D$ is the unit disc such that 
$\phi(f^p(\phi^{-1}(z)))=e^{2\pi i \alpha}z$
for some $\alpha\in \R\backslash\Q$.
In this case, $U$ is called a Siegel disc.
\item There exists an analytic homeomorphism $\phi:U\to A$  
where
$A$ is  an annulus, $A=\{z: 1<|z|<r\}, r>1$,
such that $\phi(f^p(\phi^{-1}(z)))=e^{2\pi i \alpha}z$
for some $\alpha\in \R\backslash\Q$.
In this case, $U$ is called a Herman ring.
\item There exists $z_0\in\partial U$ such that
$f^{np}(z)\to z_0$ for $z\in U$ as $n\to \infty$, but 
$f^p(z_0)$ is not
defined. In this case, $U$ is called a Baker domain.
\end{itemize}
\end{classification}
Clearly,  if $f$ is rational, then Baker domains do
not exist. If $f\in E$, then Baker domains are possible
only for $z_0=\infty$. Similarly, if $f\in P$ with pole
at $0$,
then Baker domains are possible
only for $z_0\in\{0,\infty\}$.
 
The above classification theorem is
essentially due to Cremer \cite{Cre32}
and Fatou \cite{Fat19}.
Fatou \cite[\S 56, p.\ 249]{Fat19} proved that if 
$\{f^n|_U\}$ has
only constant limit functions, then $U$ is an immediate 
attractive
basin or a Leau domain, provided
$f$ is rational. His proof shows that the
only further possibility in the case of transcendental 
functions is
that of a Baker domain. Cremer \cite[p.\ 317]{Cre32} proved 
that if
$\{f^n|_U\}$ has nonconstant limit functions, then $U$ is a 
Siegel disc
or a Herman ring. Neither Fatou nor Cremer
stated the full  classification theorem, but T\"opfer's 
remarks
\cite[p.\ 211]{Toe49} come fairly close to it.
 
We remark that when Fatou and Cremer wrote their papers, it 
was
not known yet that Siegel discs and Herman rings do actually
exist,
and they may have believed that such domains do not exist.
(Cremer \cite[p.\ 154]{Cre28} wrote that it is conjectured 
that
rational functions do not have Siegel discs but also that he
does
not see a reason for this conjecture.)
T\"opfer knew about Siegel discs, but the existence of 
Herman
rings (which T\"opfer called ``Zentrumring'') was not 
established yet.
 
In the above form the classification theorem was stated 
first
by Baker, Kotus, and L\"u 
\cite[Theorems 2.2 and 2.3]{Bak91}.
In the case of rational functions
it seems to have been given
first by Sullivan \cite{Sul82,Sul83}.
We remark that the case of an immediate
attractive basin is sometimes further
distinguished depending on whether the attracting periodic 
point
contained in it is superattracting or not. If this is the 
case,
then $U$ is called a {\em B\"ottcher domain}; otherwise, $U$
is
called a {\em Schr\"oder domain}.
The other notations are also not uniform in the literature:
Leau domains are also called {\em parabolic domains}; Herman
rings are also named after {\em Arnol'd}; and for Baker
domains the names {\em infinite Fatou component} 
\cite{Her85},
{\em essentially parabolic domain} \cite{Bak92}, and
{\em domains at} $\infty$ \cite{Dev89} are also used.
The term ``Baker domain'' seems to have been used first in
\cite{Ere90,Ere90a}.
 
Besides the papers cited already, we refer to 
\cite{Bea91,Mil90,Ste90}
for a proof of the classification theorem. Here only the 
case that
$f$ is rational is considered, but the changes necessary to 
handle
the case that $f$ is transcendental are minor.
 
We note that if $f$ is entire, then $f$ does not have Herman
rings.
In fact, a simple argument shows that $f^n\to \infty$ in
mul\-ti\-ply con\-nec\-ted
components of $F$; see \cite[p.\ 67]{Toe39}.
Baker \cite[Theorems 1 and 3]{Bak87a}
proved that analytic self-maps of the punctured plane may 
have
a Herman ring but at most one (which has period $1$).
It is not clear whether this is possible for $f\in P$.
 
\subsection{The role of the
singularities of the inverse function}  \label{inverse}
The periodic domains are  closely related to the set
of singularities of the inverse function $f^{-1}$ of $f$,
that is, the set of critical and finite asymptotic values of
$f$ and
(finite) limit
points of these values.
Denote this set by $\Si(f^{-1})$.
\begin{sing} \label{sing}
Let $f$ be a meromorphic function, and let
$C=\{U_0,U_1,\dots,U_{p-1}\}$ be a periodic cycle of 
components
of $F$.
\begin{itemize}
\item If $C$ is a cycle of
immediate attractive basins or Leau domains, then
$U_j\cap \Si(f^{-1})\neq \emptyset$ for some
$j\in\{0,1,\dots,p-1\}$.
More precisely, there exists
$j\in\{0,1,\dots,p-1\}$ such that
$U_j\cap \Si(f^{-1})\neq \emptyset$ contains a point
which is not preperiodic
or such that $U_j$ contains a periodic critical point
\rom{(}in which case $C$ is a cycle of superattractive 
basins\rom{)}.
\item If $C$ is a cycle of
Siegel discs or Herman rings, then
\!$\partial U_j\!\subset\! \overline{O^+(\Si(f^{-1})\!)}$
for all $j\in\{0,1,\dots,p-1\}$.
\end{itemize}
\end{sing}
These results were proved by Fatou \cite[\S\S 30--31]{Fat19}
for rational maps, but the proofs extend to the 
transcendental case.
 
It follows from Theorem \ref{sing} that the number of cycles
of
immediate attractive basins and Leau domains does not exceed
the number of singularities of $f^{-1}$. For transcendental
functions, $\Si(f^{-1})$ may of course be infinite (and 
simple
examples like $f(z)=e^z+z+1$ or $f(z)=e^z+z+2$
show that there may, in fact, be infinitely many cycles of
immediate attractive basins and Leau domains), but for a 
rational
function $f$ of degree $d$ there are at most $2d-2$ 
singularities
of $f^{-1}$.
 
The number of cycles of Siegel discs and Herman rings of a 
rational
function may also be bounded in terms of the degree $d$. The
sharp bound
is due to Shishikura \cite[p.\ 5]{Shi87}, who,
strengthening earlier results of Fatou
\cite[\S 30]{Fat19} and Sullivan \cite[p.\ 6]{Sul83},
proved that the number of cycles of
immediate attractive basins, Leau domains, and Siegel discs
plus twice the number of cycles of Herman rings does not 
exceed $2d-2$.
Loosely speaking,
cycles of immediate attractive basins, Leau domains, and 
Siegel discs
require one critical point, while cycles of Herman rings 
need at least
two.
(Theorem \ref{sing} gives an heuristic argument of why this 
should be
true but by no means proves it.)
A result of Shishikura's type for a class of transcendental 
entire
functions can be found in \cite[Theorem 5]{Ere90a}.
 
One may ask whether Baker domains are also related to 
singularities
of $f^{-1}$.  Examples in
\cite[Example 3]{Ere87} and \cite[p.\ 609]{Her85} show that 
a
periodic cycle of Baker domains need not contain points
of $\Si(f^{-1})$.
\begin{qu12} \label{qu12}
Let $f$ be a meromorphic function with  a cycle of Baker 
domains
that does not contain a point of $\Si(f^{-1})$. Is there 
some relation
between $\Si(f^{-1})$ and the boundaries of the domains of 
this cycle?
\end{qu12}
We will see in
Theorem \ref{asympt} in \S \ref{baker} that if 
$f^n\not\to\infty$
in a cycle of Baker domains,
then some domain in this cycle has a finite asymptotic
value on its boundary (regardless of whether there are 
singularities
of $f^{-1}$ in this cycle or not).
Question \ref{qu12} asks whether more
can be said if no points of $\Si(f^{-1})$ are in this cycle.
More specifically, one may ask the following question:
\begin{qu17} \label{qu17}
Is it possible that a meromorphic function $f$ has Baker 
domains
if $O^+(z)$ is bounded for all $z\in \Si(f^{-1})$?
\end{qu17}
 
\subsection{The connectivity of the components of the Fatou 
set}
\label{connectivity}
By definition, Siegel discs  are simply connected and
Herman rings are doubly con\-nec\-ted.
In this section we consider the connectivity of the other
components of $F$.
\begin{th5} \label{th5}
Let $f$ be a meromorphic function, and let $U$ be an 
invariant component
of $F$. Then the connectivity of $U$ has one of the values 
$1$, $2$,
or $\infty$. Here $2$ occurs only when $U$ is a Herman ring.
\end{th5}
This was proved by Fatou \cite[\S 32]{Fat19}
if $f$ is rational (see also \cite[\S 7.5]{Bea91} for this 
case)
and by Baker, Kotus, and L\"u \cite[Theorem 3.1]{Bak91a}
if $f\in M$.
Of course, the result implies that the connectivity of a 
periodic
component (of period greater than $1$) also takes one of the
values
$1$, $2$, or $\infty$ if $f$ rational.
Probably this remains true for functions in $M$ as well, but
the
proof in \cite{Bak91a}
does not seem to give this result.
\begin{qu19}
Let $f$ be a meromorphic function, and let $U$ be a periodic
component
of $F$. Is the connectivity of $U$ either $1$, $2$,
or $\infty$?
\end{qu19}
 
Baker, Kotus, and L\"u \cite[Theorem 6.1]{Bak91a} also 
proved that,
in contrast to Theorem \ref{th5}, the connectivity of
a preperiodic component may take any value if $f$ is 
rational
or if $f\in M$.
Moreover, they gave examples of functions in $M$
which have a wandering domain of any preassigned 
connectivity
\cite{Bak90b}.

For functions in $P$ and $E$, we have results stronger
than Theorem \ref{th5}.
In fact, as proved by Baker \cite[Theorem 1]{Bak87},
the connectivity of any component
of $F$ is $1$ or $2$ if $f$ is an analytic self-map of
$\C\backslash \{0\}$ and, hence, in particular if $f\in P$.
For entire functions we have the following result.
\begin{th6} \label{th6}
If $f\in E$, then any preperiodic component of $F$ is simply
connected.
\end{th6}
In other words, mul\-ti\-ply con\-nec\-ted components of $F$
are
necessarily wandering if $f\in E$.
 
Theorem \ref{th6} is an immediate consequence
of a result of Baker \cite[Theorem 1]{Bak75},
who proved that mul\-ti\-ply con\-nec\-ted components of 
$F(f)$
are bounded if $f\in E$.
In order to give a proof of Theorem \ref{th6}
(following Baker's argument), we start
with the following lemma.
With further applications in mind, this lemma is
stated in a form more general than needed for the
proof of Theorem \ref{th6}.
The results contained in it  can be found in
\cite[Lemmas 1 and  2; 23, Lemma 4.1]{Bak88}
(see also \cite[Theorem 6; 103, Proposition A.1]{Bak81}).
\begin{hypmetric} \label{hypmetric}
Let $G$ be an unbounded open set in $\C$ with at least two 
finite
boundary points, and let $g$ be analytic in $G$.
Let $D$ be a domain contained in
$G$, and suppose that $g^n(D)\subset G$
for all $n$ and that
$g^n|_D\to\infty$ as $n\to\infty$.
Then, for any compact subset $K$ of $D$, there exist
constants $C$ and $n_0$  such that
\begin{equation}
|g^n(z')|\leq |g^n(z)|^{C} \label{hyp3}
\end{equation}
for all $z,z'\in K$ and $n\geq n_0$.
If, in addition, $g(D)\subset D$, then we also have
\begin{equation}
\log \log |g^n(z)|=O(n) \label{hyp1}
\end{equation}
for all $z\in D$ as $n\to\infty$, and there exist a
constant $A>1$ and a curve
$\gamma\subset D$ tending to $\infty$
which satisfies  $g(\gamma)\subset \gamma$  such
that
\begin{equation}
|z|^{1/A}\leq |g(z)|\leq |z|^A \label{hyp2}
\end{equation}
for $z\in \gamma$.
If $\widehat{\C}\backslash G$ contains a connected set 
$\Gamma$ such that
$\{a,\infty\}\subset \Gamma$ for some $a\in \C$, then 
\rom{(\ref{hyp3}),
(\ref{hyp1})}, and \rom{(\ref{hyp2})} may be replaced by
\begin{equation}
|g^n(z')|\leq C|g^n(z)| \label{hyp6},
\end{equation}
\begin{equation}
\log |g^n(z)|=O(n) \label{hyp4},
\end{equation}
and
\begin{equation}
\frac{|z|}{A}\leq |g(z)|\leq A|z| \label{hyp5}.
\end{equation}
In particular, this is the case if $G$ is simply connected.
\end{hypmetric}
In order to prove Lemma \ref{hypmetric}, we denote by
$\Omega$ the plane punctured at two finite boundary points
of $G$. By $[z,z']_\Omega$ we denote the hyperbolic
distance of two points $z$ and $z'$ in $\Omega$.
To prove (\ref{hyp3}), we note that
\[
[g^n(z),g^n(z')]_\Omega\leq
[g^n(z),g^n(z')]_G\leq [g^n(z),g^n(z)]_{g^n(D)}\leq [z,z']_D
\]
and that the hyperbolic metric $\rho_\Omega(z)$ satisfies
\[
\rho_\Omega(z)\sim\frac{c}{|z|\log|z|}
\]
for some positive constant $c$ as $|z|\to\infty$;
see \cite[\S 1.8]{Ahl73}.
It follows that if
$|g^n(z')|\geq |g^n(z)|$, then
\[
[g^n(z),g^n(z')]_\Omega
\geq
\frac{c}{2}
\int_{|g^n(z)|}^{|g^n(z')|}\frac{dt}{t\log t}
=
\frac{c}{2}
\log \left( \frac{\log|g^n(z')|}{\log|g^n(z)|} \right)
\]
for sufficiently large $n$ so that (\ref{hyp3})
holds with
\[
C = \exp \left( \frac{2}{c}
\max_{z,z'\in K} [z,z']_D \right).
\]
In order to prove (\ref{hyp2}),
we choose $\sigma$ as a curve in $D$ that connects a point
$z_0\in D$ with $g(z_0)$ and define
$\gamma = \bigcup_{n=0}^\infty g^n(\sigma)$.
Then (\ref{hyp2}) can be deduced from (\ref{hyp3}) if we 
choose
$K=\sigma\cup g(\sigma)$. Similarly, choosing $z'=g(z)$ in 
(\ref{hyp3}),
we have
\[
|g^n(z)|=|g^{n-1}(z')|\leq |g^{n-1}(z)|^C
\leq |g^{n-2}(z)|^{C^2}\leq\dotsb
\]
for large $n$, and (\ref{hyp1}) follows
by induction.
 
To prove (5)--(7),
we proceed as above but define $\Omega$ as the component of
$\widehat{\C}\backslash \Gamma$ that contains $G$.
Then $\rho_\Omega(z)\geq c/|z|$ for some positive constant 
$c$
and all sufficiently large $z\in \Omega$. The arguments
used above  to prove
(2)--(4)
now yield (5)--(7).
This completes the proof of Lemma \ref{hypmetric}.
 
We now  prove Theorem \ref{th6}.
Suppose that $U$ is a mul\-ti\-ply con\-nec\-ted
component of $F$, and let $\sigma\subset U$ be a curve that 
is not
null-homotopic in $U$.
Define $\sigma_n=f^n(\sigma)$. Then $\sigma_n$ is not 
null-homotopic
in $U_n$ so that $U_n$ is mul\-ti\-ply con\-nec\-ted for all
$n$.
It is not difficult to see that $f^n|_U\to \infty$
as $n\to\infty$.
By (\ref{hyp3}), there exists a sequence $(r_n)$
tending to $\infty$ and a constant $C$ such that
$\sigma_n\subset \Ann(r_n,r_n^C)$, where
$\Ann(r,R)=\{z:r<|z|<R\}$ if $0\leq r<R$.
By Theorem \ref{rep} there exists a periodic point $z_0$
contained in $\Int(\sigma)$, the interior of
$\sigma$. It follows that for all $n$,
some point of the periodic
cycle to which $z_0$ belongs is contained in 
$\Int(\sigma_n)$.
Hence $D(0,r_n)\subset \Int(\sigma_n)$
for sufficiently large $n$, where $D(a,r)=\{z:|z-a|<r\}$ for
$r>0$ and $a\in\C$.
 
Suppose now that $U$ is preperiodic. Replacing $U$ by $U_m$
and $f$ by $f^n$ for suitable values of $m$ and $n$, we
may assume without loss of generality that $U$ is invariant.
We deduce from Lemma \ref{hypmetric} that
there exist a constant $A$ and
a curve $\gamma$  tending to $\infty$ such
that $|f(z)|\leq |z|^A$ for $z\in \gamma$.
For sufficiently large $n$, $\sigma_n$ intersects $\gamma$;
that is,  we can find $w_n\in\sigma_n\cap\gamma$.
Denoting by $M(r,f)$ the maximum modulus of $f$, that is,
$M(r,f)=\max_{|z|=r}|f(z)|$,
we deduce that
\[
M(r_n,f)
\leq
\max_{z\in\sigma_n}|f(z)|
\leq
r_{n+1}^C
\leq
\left(\min_{z\in\sigma_n}|f(z)|\right)^C
\leq
|f(w_n)|^C
\leq
|w_n|^{CA}
\leq
(r_n)^{C^2A}.
\]
This is a contradiction to the hypothesis that $f$ is 
transcendental
and thus completes the proof of Theorem \ref{th6}.
 
The argument used above actually shows that certain classes
of entire functions do not have mul\-ti\-ply con\-nec\-ted 
domains
of normality at all.
For example, using this method, one can obtain the following
result.
\begin{th7} \label{th7}
Suppose that $f\in E$ and that for all $\varepsilon>0$ there
exists
a curve $\gamma$ tending to $\infty$ such that
$|f(z)|\leq M(|z|^\varepsilon,f)$ for $z\in \gamma$.
Then all components of $F$ are simply connected.
In particular, this is the case  if
$\log |f(z)|=O(\log |z|)$ as $z\to\infty$ through some path.
\end{th7}
On the other hand, examples of entire functions with
mul\-ti\-ply con\-nec\-ted
components of the Fatou set are known. The first example was
constructed in \cite{Bak63}; further examples can be found 
in
\cite{Bak76,Bak85,Hin92a}.
 
Baker \cite{Bak85} gave an example of a transcendental 
entire
function with an infinitely connected (and hence wandering)
domain of normality. In the other examples cited, it is not 
clear
what the connectivity of the
mul\-ti\-ply con\-nec\-ted
components is.
 
\begin{qu3}
Is there
an entire transcendental
function whose Fatou set  has
mul\-ti\-ply con\-nec\-ted components of
finite connectivity?
\end{qu3}
 
\subsection{Wandering domains}  \label{wander}
The first example of an entire function with a wandering 
domain
was given by Baker \cite{Bak76}. We remark that the 
existence
of wandering domains follows directly from his results
in \cite{Bak63} (existence of
mul\-ti\-ply con\-nec\-ted domains of normality)
and \cite{Bak75}
(mul\-ti\-ply con\-nec\-ted domains of normality are 
bounded),
but \cite{Bak76} was written before \cite{Bak75} (although
published later).
 
Since then, many other examples have been constructed.
The elementary examples
\[
f_1(z)=z-1+e^{-z}+2\pi i
\]
and
\[
f_2(z)=z+\lambda \sin(2\pi z)+1
\]
where $1+2\pi\lambda=e^{2\pi i\alpha}$ for suitable real
numbers $\alpha$ have been given by Herman (\cite[p.\ 106;
130, p.\ 414]{Her84};
see also \cite[pp.\ 564, 567]{Bak84}).
While the wandering domain in Baker's example \cite{Bak76}
is mul\-ti\-ply con\-nec\-ted, $f_1$ and $f_2$ have simply 
connected
wandering domains.
 
To prove that $f_1$ has a wandering domain, define
$g(z)=z-1+e^{-z}$. (The function $g$ arises if we apply
Newton's method to $h(z)=e^z-1$, that is, 
$g(z)=z-h(z)/h'(z)$.)
For $k\in \Z$ we define $z_k=2\pi k i$. Then $z_k$ is a
superattracting fixed point of $g$. Denote by $U_k$ the
immediate attractive basin of $z_k$, that is, the
component of $F(g)$ that contains $z_k$.
One can show that $J(g)=J(f_1)$; compare
\cite[Lemma 4.5]{Bak84}.
It follows that $U_k$ is also a component of $F(f_1)$,
and we clearly have $f(U_k)\subset U_{k+1}$; that is,
$U_k$ is wandering. It is not difficult to see that
$U_k$ is simply connected for all $k$.
 
The proof that $f_2$ has a simply connected
wandering domain is similar.
Here $\alpha$ is chosen such that $z+\lambda\sin(2\pi z)$
has a Siegel disc at zero. In this example, we obtain
wandering domains $U_k$ containing $k\in \Z$.
Here the $U_k$ have the additional feature that
$f|_{U_k}$ is univalent.
Different examples of wandering domains with
this property have been constructed by Eremenko
and Lyubich \cite[Example 2]{Ere87}.
 
An example similar to $f_1$ and $f_2$ is given by
\[
f_3(z)=z+\lambda \sin z
\]
where $\lambda\in \R$ is chosen so that the forward orbit
of each critical point consists only of critical points.
For a discussion of this example, see
\cite[p.\ 222; 56, p.\ 290; 60, p.\ 52]{Dev86}.
Other examples of wandering domains with various
additional properties have also been given.
For example, Baker \cite[Theorem 1]{Bak85} (see also
\cite[Theorem 5.2]{Bak84}) has shown that the order
of an entire function with wandering domains may take
any value.
 
In all examples mentioned so far, the iterates tend to
$\infty$ in the wandering domain.
It is well known 
(see \cite[Lemma 2.1; 48, p.\ 317; 71, \S28]{Bak91a})
that there cannot exist nonconstant limit
functions of $\{f^n|_U\}$ if $U$ is a wandering domain of
a meromorphic function $f$.
Eremenko and Lyubich \cite[Example 1]{Ere87} have
constructed an entire function $f$ with a wandering
domain $U$ such that the set of limit functions of
$\{f^n|_U\}$ contains an infinite number of
finite constants.
In this example,
the constant limit functions  have $\infty$ as a limit
point; that is, $\infty$ is also a limit function of
$\{f^n|_U\}$. It is a well-known open problem
\cite[Problems 2.77 and 2.87]{Bra89} whether
this is always the case.
\begin{qu4}
Let  $U$ be a wandering domain of the transcendental
meromorphic function $f$. Does there exist a sequence
$(n_k)$ such that $f^{n_k}|_U\to\infty$ as $k\to\infty$?
\end{qu4}
We remark that it has been shown in \cite{Ber92b}
that if $f\in E$ and $U$ is a wandering
domain of $f$, then all finite limit functions
of $\{f^n|_U\}$ are contained in the derived set
of $O^+(\Si(f^{-1}))$.
 
Finally, we mention that Baker, Kotus, and L\"u
\cite[\S 6]{Bak90b}
have modified the method of Eremenko and Lyubich
to construct a
function  $f\in M$ which has a mul\-ti\-ply con\-nec\-ted
wandering domain $U$ of preassigned connectivity
such that
the limit set of  $\{f^n|_U\}$ contains infinitely many
finite constants.
 
\subsection{Classes of functions without wandering domains}
\label{nowander}
One of the most important results in the iteration theory of
rational
functions is the following theorem of Sullivan 
\cite{Sul82,Sul85}.
\begin{sullivan} \label{sullivan}
Rational functions do not have wandering domains.
\end{sullivan}
Together with Theorems \ref{classification} and \ref{sing},
this leads to a fairly complete description of the
iterative behavior of rational functions on the Fatou set.
 
Sullivan's theorem has been extended to various classes
of transcendental functions.
We mention the following classes:
\begin{itemize}
\item $S=\{f:f$ has only finitely many critical and 
asymptotic values$\}$;
\item $F=\{f:f$ has a representation of the form 
$f(z)=z+r(z)e^{p(z)}$\newline
\qquad\qquad\quad where
 $r$ is rational and $p$ is a polynomial$\}$;
\item $N=\{f:f$ has finite order and 
$f'(z)=r(z)e^{p(z)}(f(z)-z)$\newline
\qquad\qquad\quad
where $r$ is rational and $p$ is a polynomial$\}$;
\item $R=\{f:
f'(z)=r(z)(f(z)-z)^2$ or $f'(z)=r(z)(f(z)-z)(f(z)-\tau)$
\newline
\hfill where $r$ is rational and $\tau\in\C\}$.
\end{itemize}
The names of the different classes are somewhat arbitrary. 
According
to Eremenko and Lyubich \cite[p.\ 624]{Ere90} $S$ was
chosen in honor of {\em S}peiser, who introduced
this class in a different context.
Class $N$ is of interest in connection with {\em N}ewton's 
method
(compare \S \ref{newton}),
and class $R$ consists of solutions of certain {\em R}iccati
equations
(but is also of interest for Newton's method).
\begin{nowd} \label{nowd}
Functions in $S$, $F$, $N$, and $R$ do not have wandering 
domains.
\end{nowd}
We note that all these classes contain the class
of rational functions so that
Theorem \ref{nowd} may be considered as a generalization of
Theorem \ref{sullivan}.
The result that meromorphic functions in $S$ do not have 
wandering
domains was proved by Baker, Kotus, and L\"u \cite{Bak92}.
This result had been obtained earlier
by Eremenko and Lyubich \cite{Ere84,Ere90a}
and Goldberg and Keen \cite{Gol86} for $S\cap E$
and by Keen \cite{Kee88}, Kotus \cite{Kot87}, and Makienko 
\cite{Mak88}
for $S\cap P$ (and, in fact, for the corresponding class of
analytic self-maps of the punctured
plane).
For other subclasses of $S$, this had been proved by
Baker \cite[Theorem 6.2]{Bak84}
and Devaney and Keen \cite[p.\ 72]{Dev89}.
 
The result that functions in $F$ do not have wandering
domains was proved by Stallard \cite{Sta91}.
The result for the classes $N$ and $R$ can be found in
\cite{Ber92a} and \cite{Ber92c}, respectively.
In \cite{Ber92c} the nonexistence of wandering domains
is also proved for solutions of certain other differential
equations.
 
The proofs in the papers cited above depend
crucially on the fact that if $f$ is in one of the above 
classes,
then there exist only finitely many singularities of 
$f^{-1}$
that are not contained in preperiodic components.
(This is clear for $f\in S$ and easy to see for
$f\in N$ and $f\in R$, but the proof for $f\in F$ is more
involved; see \cite{Sta91}.)
It follows that if $U$ is a wandering domain, then there
exists $n_0$ such that $U_n\cap \Si(f^{-1})=\emptyset$
for $n\geq n_0$. Now two cases have to be distinguished:
\newcounter{counter}
\begin{list}
{(\roman{counter})}{\usecounter{counter}}
\item $U_n$ is simply connected for all $n\geq n_0$.
\item $U_m$ is
mul\-ti\-ply con\-nec\-ted for some $m\geq n_0$.
\end{list}
In case (i), one uses the
ideas of Sullivan \cite{Sul85}, who introduced
quasiconformal mappings into the subject. We sketch the
argument very briefly.
Consider a quasiconformal homeomorphism of $U_{n_0}$ with
complex dilatation $\mu$.
Then $\mu$ can be extended to $\widehat{\C}$ in such a way 
that
$\mu(f(z))=\mu(z)f'(z)/\overline{f'(z)}$ for all $z\in \C$.
Then there exists a quasiconformal homeomorphism $\Phi$ of
$\widehat{\C}$ that fixes $0$, $1$, and $\infty$ and whose
complex dilatation is $\mu$.
We now consider $f_\Phi=\Phi\circ f\circ \Phi^{-1}$
and observe that if $f$ is in one of the classes under
consideration, then so is $f_\Phi$. This sharply limits
the possibilities for $f_\Phi$, and---loosely speaking---a 
contradiction is obtained from
the fact that there are many
quasiconformal homeomorphisms of $U_{n_0}$ and
hence many functions $\mu$
but not so many functions $f_\Phi$.
For the details we refer to \cite{Bak84,Bak92,Bea91}.
 
In  case (ii) it is not difficult to  obtain a contradiction
to Theorem \ref{th7} if $f$ is entire and contained in
$S$, $F$, or $N$.
A result similar to Theorem \ref{th7} can still be obtained
if $f$ has finitely many poles, and this has been used to
rule out case (ii) for meromorphic functions in $F$ and $N$;
see \cite{Ber92a} and \cite{Sta91} for details.
For meromorphic functions in $S$, a different but still
fairly elementary argument has been used; see \cite{Bak92}.
The proof that case (ii) cannot occur for $f\in R$ is less
elementary but uses results of Shishikura \cite{Shi90}
obtained by quasiconformal surgery. We refer to 
\cite{Ber92c}
for the details.
 
Besides the classes contained in Theorem \ref{nowd}, there
are some other classes of functions known to have no
wandering domains.
We mention that if $g$ is an analytic self-map of
$\C\backslash\{0\}$, then there exist entire functions $f$
satisfying $\exp\circ f=g\circ \exp$. Thus results obtained
for
analytic self-maps of $\C\backslash\{0\}$ may be used to 
obtain
results for entire functions $f$ that admit a representation
of the above form. For example,  one can prove using these
ideas that if $p$ and $q$ are
polynomials, then $f(z)=p(e^z)+q(e^{-z})$ does not have 
wandering
domains. For this and related results
we refer to \cite{Bak84,Kee88,Kot87,Mak88}.
 
If we combine the already mentioned result in \cite{Ber92b} 
that
all finite limit functions of $\{f^n|_U\}$
are contained in the derived set of $O^+(\Si(f^{-1}))$ if
$f\in E$ and if $U$ is a wandering domain of $f$
with Theorem \ref{ere90} in \S \ref{nobakdom},
we also obtain
some classes of entire functions without
wandering domains.
We note that this is a fairly elementary way to obtain
classes of functions without wandering domains, while
the proofs of Theorems \ref{sullivan} and \ref{nowd}
use quasiconformal mappings.
\begin{qu6}
Is there a proof of Theorem \ref{sullivan} (and
Theorem \ref{nowd}) that does not use
quasiconformal mappings?
\end{qu6}
Some of the
results concerning the nonexistence of wandering
domains suggest that there are relations between wandering 
domains
and singularities of the inverse function.
In fact, similarly to \S
\ref{inverse}
where we said that
periodic components of the Fatou set
(seem to) require one respectively two
singularities, one may ask whether wandering domains require
infinitely many of them, in a sense which still has to be 
made precise.
More specifically, one may ask the following questions.
\begin{qu5} \label{qu5}
Can a meromorphic function $f$ have wandering domains
if all (or all but finitely many) points
of $\Si(f^{-1})$ are contained
in preperiodic domains?
\end{qu5}
\begin{qu13}
Let $f$ be a meromorphic function with a wandering domain
$U$ such that $U_n\cap\Si(f^{-1})=\emptyset$ for all 
$n\geq 0$.
Is there some relation between $\partial U_n$ and
$\Si(f^{-1})$?
\end{qu13}
 
\subsection{Baker domains} \label{baker}
The first example of an entire function with
a Baker domain was already given by Fatou 
\cite[Example I]{Fat19},
who considered the function
\[
f(z)=z+1+e^{-z}
\]
and proved that $f^n(z)\to\infty$ as $n\to\infty$ for 
$\re z>0$,
that is, the right half-plane is contained in an invariant 
Baker domain.
An example of a Baker domain of higher period was given by 
Baker,
Kotus, and L\"u \cite[p.\ 606]{Bak91a}, who showed that
the function $f(z)=1/z-e^z$ has a cycle $\{U_0,U_1\}$ of 
Baker
domains such that $f^{2n}|_{U_0}\to \infty$ and
$f^{2n}|_{U_1}\to 0$ as $n\to\infty$.
 
We
list some general properties of Baker domains.
Let $\{U_0,U_1,\dots,U_{p-1}\}$ be a periodic cycle of Baker
domains, and denote by $z_j$ the limit corresponding to
$U_j$, that is, $f^{np}(z)\to z_j$ for $z\in U_j$ as 
$n\to\infty$.
Clearly, $f(z_j)=z_{j+1}$ if $z_j\neq\infty$.
(Here, by definition, $z_p=z_0$.)
It follows that there exists at least one 
$j\in\{0,1,\dots,p-1\}$
such that $z_j=\infty$, and for all
$j\in\{0,1,\dots,p-1\}$ there exists
$l=l(j) \in\{0,1,\dots,p-1\}$ such that
$f^l(z_j)=\infty$.
 
The $U_j$ contain curves $\gamma_j$ tending to $z_j$ such 
that
$f^p(\gamma_j)\subset\gamma_j$
and $f^p(z)\to z_j$ as $z\to z_j$ in $\gamma_j$.
To see this, we proceed as  in the proof of Lemma 
\ref{hypmetric}
and
choose $w_0\in U_0$ and a curve  $\sigma\subset U_0$
that joins $w_0$ and
$f^p(w_0)$. We
define $\gamma_0=\bigcup_{n=0}^\infty f^{np}(\sigma)$
and $\gamma_j=f^j(\gamma_0)$ for $j\in\{1,2,\dots,p-1\}$.
Then the $\gamma_j$ have the desired properties. Moreover,
$f(z)\to z_{j+1}$ as $z\to z_j$ in $\gamma_j$.
We deduce that if $z_j=\infty$, then $z_{j+1}$ is an
asymptotic value of $f$, the asymptotic path being contained
in $U_j$.
 
We collect some of the above observations
in the following theorem.
\begin{asympt} \label{asympt}
Let $f$ be a meromorphic function, and
let $\{U_0,U_1,\dots,U_{p-1}\}$ be a periodic cycle of Baker
domains of $f$. Denote by $z_j$ the limit corresponding to
$U_j$, and define $z_p=z_0$.
 Then $z_j\in \bigcup_{n=0}^{p-1}f^{-n}(\infty)$ for
all $j\in\{0,1,\dots,p-1\}$, and $z_j=\infty$ for at least 
one
$j\in\{0,1,\dots,p-1\}$. If $z_j=\infty$, then $z_{j+1}$ is 
an
asymptotic value of $f$.
\end{asympt}
\begin{co1} \label{co1}
If $f$ has a cycle
$\{U_0,U_1,\dots,U_{p-1}\}$ of Baker domains such that
$f^n|_{U_0}\to\infty$, then
$\infty$ is an asymptotic value of $f$.
In particular, this is the case
if $f$ has an invariant Baker domain.
\end {co1}
\begin{co2} \label{co2}
If $f$ has a cycle
$\{U_0,U_1,\dots,U_{p-1}\}$ of Baker domains such that
$f^n|_{U_0}\not\to\infty$, then
$f$ has a finite asymptotic value.
\end{co2}
Corollary  \ref{co1} can be found in \cite[p.\ 75]{Dev89}
for maps with polynomial Schwarzian derivative.
 
Lemma \ref{hypmetric} gives additional information about the
asymptotic paths $\gamma_j$  and also answers the question 
how
fast $f^{np}(z)$ approaches
$z_j$ for $z\in U_j$.
In fact,
if $U_j$, $z_j$, and $\gamma_j$ are as above,
then $|z|^{1/A}\leq|f^p(z)|\leq |z|^A$
for $z\in \gamma_j$
and
$\log\log |f^{pn}(z)|=O(n)$ for $z\in D$
if $z_j=\infty$.
If $U_j$ is simply connected and $z_j=\infty$, then
we even have
$|z|/A\leq |f^p(z)|\leq A|z|$
for $z\in\gamma_j$ and
$\log |f^{pn}(z)|=O(n)$ for $z\in D$.
Similar results may be obtained if $z_j\neq\infty$.
 
As already mentioned after Theorem \ref{sing},
periodic cycles of Baker domains need not
contain a singularity of $f^{-1}$.
However, we have the following result.
\begin{ber92b} \label{ber92b}
If $f\in N$ or $f\in F$, then any periodic cycle of Baker 
domains
contains a point  of $\Si(f^{-1})$.
\end{ber92b}
This result was proved in \cite{Ber92a} for $f\in N$, but 
the
proof extends to the case that $f\in F$.
The proof of Theorem \ref{ber92b}
is fairly analogous to the proof that
functions in $N$ (and $F$) do not have wandering domains.
Therefore, it seems likely that
the conclusion of Theorem \ref{ber92b}
remains valid for functions in $R$.
(This is certainly so for cycles of simply connected
Baker domains, but in the mul\-ti\-ply con\-nec\-ted case 
some
modification of the argument will have to be made.)
For functions in $S$ we have a stronger result; see 
Corollary \ref{co4}
in \S \ref{nobakdom}.
 
One way to prove that cycles of Leau domains contain a 
singularity of
$f^{-1}$ is based on the solution of Abel's functional 
equation
(cf. \cite[\S 7]{Mil90}).
Hinkkanen \cite[Theorem 2]{Hin92}
has shown that in certain cases
this argument may  also be used to prove that
Baker domains contain singularities of $f^{-1}$.
 
\subsection{Classes of functions without Baker domains}  
\label{nobakdom}
Eremenko and Lyubich \cite{Ere90a}
considered the class
\[
B=\{f:\Si(f^{-1}) \text{\ is bounded}\}
\]
and proved the following result.
\begin{ere90} \label{ere90}
If $f\in E\cap B$, then there does not exist a component
$U$ of $F(f)$ such that $f^n|_U\to\infty$ as $n\to\infty$.
\end{ere90}
\begin{co3} \label{co3}
If $f\in E\cap B$, then $f$ does not have Baker domains.
\end{co3}
We note that the conclusion of Corollary \ref{co3} does not 
hold
in general
for $f\in M\cap B$.
As an example,  consider $f(z)=1/z-e^{z}$.
As already mentioned above,
Baker, Kotus, and L\"u \cite[p.\ 606]{Bak91a}
proved that
$f$ has a Baker domain of period $2$,
and it is easy to check that $f\in M\cap B$.
In this example, the critical values of $f$ accumulate at 
$0$,
which is also one of the limits corresponding to the cycle 
of
Baker domains.
 
The following result is a generalization of Corollary 
\ref{co3}
to meromorphic functions.
\begin{bakdom} \label{bakdom}
Let $f$ be a meromorphic function, and
let $\{U_0,U_1,\dots,U_{p-1}\}$ be a periodic cycle of Baker
domains of $f$.
Then $\infty$ is in the derived set of
\[
\bigcup_{j=0}^{p-1} f^j(\Si(f^{-1})).
\]
\end{bakdom}
\begin{co4} \label{co4}
Functions in $S$ do not have Baker domains.
\end{co4}
Combining Corollary \ref{co4}  with Theorem \ref{nowd},
we see that the iteration of functions in $S$ is in many 
ways
analogous to that of rational functions and may thus be 
analyzed in
a similar way.
 
For example, these results allow us to prove that the 
functions
$\lambda z e^z$, $\lambda e^z/z$, and $\lambda\tan z$
satisfy  $J=\widehat{\C}$
for certain  values of $\lambda$, as mentioned in
\S \ref{elem}.
In fact, all these functions are in $S$ and hence do
not have wandering or Baker domains by Theorem \ref{nowd}
and Corollary \ref{co4}. For suitably chosen values of
$\lambda$ we can achieve that the points of $\Si(f^{-1})$
are either  contained in $O^-(\infty)$, or they are
preperiodic but not periodic.
In view of Theorem \ref{sing} this implies that
there are no immediate attractive basins, no Leau domains,
no Siegel discs, and no Herman rings. Hence $J=\widehat{\C}$
for these $\lambda$.
The above argument also shows that $J(e^z)=\widehat{\C}$.
 
For the proof of Theorem \ref{bakdom} we need the following 
lemma.
\begin{classB}  \label{classB}
Suppose $f\in B$, $p\geq 1$, and $0\notin O^-(\infty)$.
Then there exist a positive constant $R$ and a curve 
$\Gamma$
connecting $0$ and $\infty$ such that $|f^p(z)|\leq R$ for
$z\in \Gamma$.
\end{classB}
We show first that if $r$ is sufficiently
large, then there exists a curve $\Gamma$ connecting 
$\infty$
with some point in $\C$ such that $|f(z)|=r$ for 
$z\in\Gamma$.
In fact, otherwise the components of $f^{-1}(D(0,r))$
are bounded for arbitrarily large $r$. Hence we can find
$r_1$ and $r_2$ satisfying 
$0<r_1<r_2$ and $\Si(f^{-1})\subset D(0,r_1)$
such that $f^{-1}(D(0,r_2))$ has a bounded component which 
contains
at least two components of
$f^{-1}(D(0,r_1))$. It follows that
$f^{-1}(D(0,r_2))$ contains a component of
$f^{-1}(\Ann(r_1,r_2))$, which is at least triply connected.
By the Riemann-Hurwitz formula
(see, e.g., \cite[\S 5.4]{Bea91}), this component contains a
critical point of $f$; that is, $\Ann(r_1,r_2)$ contains
a critical value of $f$, contradicting the choice of $r_1$.
Now we choose $r$ sufficiently large
and a corresponding curve $\Gamma$
such that $\partial D(0,r)\cap O^-(\infty) =\emptyset$. Then
there exists $R>0$ such that
$|f^{p-1}(z)|\leq R$ for $|z|=r$. We deduce
that $|f^p(z)|\leq R$ for $z\in\Gamma$.
Increasing $R$ if necessary, we may assume that $\Gamma$
connects $0$ and $\infty$.
This completes the proof of Lemma \ref{classB}.
 
To prove Theorem \ref{bakdom}, we assume without loss of
generality that
$f^{np}|_{U_0}\to\infty$ as $n\to\infty$.
Suppose that the conclusion
of the theorem is false; that is, there exists a punctured
neighborhood $N_0$ of $\infty$ which does
not contain points of
$\bigcup_{j=0}^{p-1} f^j(\Si(f^{-1}))$.
In particular, this implies that $f\in B$.
We may assume without loss of generality that 
$0\notin O^-(\infty)$
so that the hypotheses of Lemma \ref{classB} are satisfied.
With $R$ and $\Gamma$ as in the conclusion of Lemma 
\ref{classB},
we may suppose that $N_0=\{z:|z|>R\}$.
In addition, we choose $R>|f(0)|$.
 
Suppose now that $w_0\in U_0\cap N_0$, and define 
$w_1=f^p(w_0)$.
We may assume that $w_1\in U_0\cap N_0$, because otherwise
we may replace $w_0$ by $f^{pn}(w_0)$ for a sufficiently
large $n$. We introduce the abbreviation $g=f^p$.
If $R$ has been chosen large enough, then the branch
of $g^{-1}$ satisfying $g^{-1}(w_1)=w_0$ may be
continued analytically in $N_0$.
We define $u_0=\log w_0$ and $u_1=\log w_1$ for arbitrary
branches of the logarithm.
Then $\Phi=\log\circ g^{-1}\circ \exp$ may be defined as a
single-valued function in the half-plane 
$H=\{z:\re z>\log R\}$
such that $\Phi(u_1)=u_0$.
Because $|g(z)|\leq R$ for $z\in \Gamma$, we have
$\Phi(u)\notin\log\Gamma$
for all $u\in H$ and any branch of the logarithm.
Hence $\Phi(H)$ does not contain a disc of radius greater
than $\pi$ so that
\begin{equation}
|\Phi'(u)|\leq \frac{\pi}{B(\re u -\log R)} \label{Phi'}
\end{equation}
where $B$ is Bloch's constant. (We do
not need any estimate for $B$ here, just Bloch's theorem
that $B>0$. Instead, we could also work with Landau's 
constant.)
In terms of $g$, we find that
\[
|g'(w)|\geq \frac{B|g(w)|(\log|g(w)|-\log R)}{\pi|w|}.
\]
Now we define $w_n=g^n(w_0)$ for $n\geq 2$, and we may 
assume
that $w_n\in U_0\cap N_0$ for all $n$.
(Otherwise,
we may replace $w_0$ again by $g^m(w_0)$ for a sufficiently
large $m$.) Then
\[
|g'(w_n)|\geq c\frac{|g(w_n)|\log|g(w_n)|}{|w_n|}
=c\frac{|w_{n+1}|\log |w_{n+1}|}{|w_n|}
\]
for some positive constant $c$ and all $n\geq 1$.
Hence
\[
|(g^n)'(w_0)|
=
\left|\prod_{j=0}^{n-1}g'(w_j)\right|
\geq
\prod_{j=0}^{n-1}
c\frac{|w_{j+1}|\log |w_{j+1}|}{|w_j|}
=
\frac{|w_n|}{|w_0|}
\prod_{j=0}^{n-1}
c\log |w_{j+1}|
\]
so that
\begin{equation}
\frac{|(g^n)'(w_0)|}{|w_n|}\to \infty    \label{g^n'}
\end{equation}
as $n\to\infty$.
We may assume that there exists a region $\Omega$ containing
$w_0$ and $w_1$
such that $g^n(\overline{\Omega})\subset U_0\cap N_0$ for 
all $n$.
We apply Lemma \ref{hypmetric} for 
$G=D=\bigcup_{n=0}^\infty g^n(\Omega)$
and note that $\Gamma$ satisfies the hypotheses of this 
lemma.
Hence 
$g^n(\Omega)\subset D(0,C |w_n|)\subset D(w_n,(C+1)|w_n|)$
for some constant $C$ by
(\ref{hyp6}). It follows that if we choose $r>0$ such that 
the disc
around $w_0$ of radius $r$ is contained in $\Omega$, then
\[
|(g^n)'(w_0)|\leq \frac{(C+1)|w_n|}{r}
\]
by Schwarz's lemma.
This contradicts (\ref{g^n'}) and completes
the proof of Theorem \ref{bakdom}.
 
The above proof uses some of the ideas introduced by 
Eremenko and Lyubich
\cite{Ere90a} to prove Theorem \ref{ere90}.
We sketch their proof of Theorem \ref{ere90}.
First, we define again $N_0=\{z:|z|>R\}$.
Using the methods of Lemma \ref{classB},
one can show that the components of $f^{-1}(N_0)$ are simply
connected
and unbounded if $R$ is sufficiently large.
This implies that if $D$ is a component of $f^{-1}(N_0)$
and if $A=\log D$ for some branch of the logarithm,
then $\Psi=\log\circ f\circ\exp$ is a conformal map from
$A$ onto $\Psi(A)$.
We define $\Phi=\Psi^{-1}$ and find again that
(\ref{Phi'}) holds.
(Here we may replace Bloch's constant by Koebe's constant,
which is equal to $\frac{1}{4}$, because $\Phi$ is 
univalent.)
Suppose now that $w_0\in F$ and $f^n(w_0)\to\infty$.
Define $u_0=\log w_0$.
Similar to (\ref{g^n'}), we find that
$|(\Psi^n)'(u_0)|\to\infty$.
On the other hand, if $U$ is a sufficiently small 
neighborhood
of $u_0$, then $\Psi^n(U)$ cannot contain a disc of radius 
larger
than $\pi$. This is a contradiction to Koebe's (or Bloch's) 
theorem.
 
\subsection{Completely invariant domains}    \label{compinv}
Recall that a set $S$ is called completely invariant
with respect to the meromorphic function $f$ if 
$O(S)\subset S$.
One may ask in which cases a component of
the Fatou set of a meromorphic function
can be completely invariant. Such components
are also called completely invariant domains.
 
It is classical that a rational function has at most
two completely invariant domains 
\cite[Theorem 5.6.1]{Bea91}.
Here the number two
is best possible, as shown by the simple example $f(z)=z^2$.
 
For transcendental entire functions, we have the following 
result of Baker
\cite{Bak70a}.
\begin{compinv1} \label{compinv1}
If $f\in E$, then $f$ has at most one completely invariant 
domain.
\end{compinv1}
It is easy to find transcendental entire functions which
have a completely
invariant domain, for example, $f(z)=\lambda e^z$ has
this property if $0<\lambda<1/e$.
 
We mention the following question of Baker.
\begin{qu14}
Suppose $f\in E$ has a
completely invariant domain $U$. Do we have $F(f)=U$?
\end{qu14}
Some results supporting the conjecture that the answer is 
``yes''
can be found in \cite[Theorem 2]{Bak75} and
\cite[\S 6]{Ere90a}.
In particular, it is shown in \cite[Theorem 6]{Ere90a}
that this is the case if $f\in S\cap E$.
We also note that Theorem \ref{compinv1}
can be deduced from \cite[Theorem 2; 70, Lemma 11]{Bak75}.
 
Less is known about
completely invariant domains of meromorphic functions.
The example $f(z)=\tan z$ where $J=\R\cup\{\infty\}$ and 
where
the upper and lower half-plane
are completely invariant shows that there may be two
completely invariant domains.
\begin{qu15} \label{qu15}
Let $f$ be a meromorphic function. Can $f$ have more than
two completely invariant domains?
\end{qu15}
A partial result was obtained by Baker, Kotus,
and L\"u \cite[Theorem 4.5]{Bak91a}.
\begin{compinv2} \label{compinv2}
If $f\in S$, then $f$ has at most two completely invariant 
domains.
\end{compinv2}
If the answer to Question \ref{qu15} is ``no'', one may also
ask
the following question.
\begin{qu16}
Suppose a meromorphic function $f$ has two
completely invariant domains $U_1$ and $U_2$. Do we have 
$F(f)=U_1\cup U_2$?
\end{qu16}
 
\section{Properties of the Julia set} \label{juliasets}
\subsection{Cantor sets and real Julia sets}  \label{cantor}
For rational functions the Julia set is often a Cantor set.
(By definition a closed subset of $\widehat{\C}$ is called a
{\em Cantor set} if it is perfect and totally disconnected.)
For example, if $\lambda$ is not contained in the Mandelbrot
set,
then $J(z^2+\lambda)$ is a Cantor set.
 
For $f\in M$ it is also possible that $J(f)$ is a Cantor 
set.
In fact, it was shown by Devaney and Keen 
\cite[p.\ 62]{Dev89}
that this is the case for $f(z)=\lambda\tan z$ if
$-1<\lambda<1$ and $\lambda\neq 0$.
 
The following result, which
is an immediate consequence of Theorem \ref{th6}, says
that this cannot happen for transcendental entire functions
\cite[p.\ 278, Corollary]{Bak75}.
\begin{continua}
If $f\in E$, then $J(f)$ contains nondegenerate continua.
\end{continua}
For rational functions it is possible that the Julia set
is a circle or a straight line; for example,
$J(\frac{1}{2}(z-1/z))=\R\cup\{\infty\}$.
This may also happen for $f\in M$. In fact, we have
$J(\lambda \tan z)=\R\cup\{\infty\}$ if $\lambda\geq 1$;
see \cite[pp.\ 60--61]{Dev89}.
For  $f\in E$ this is impossible, as shown by the following
result of T\"opfer \cite[\S 3]{Toe39}.
\begin{isoarc} \label{isoarc}
If $f\in E$, then $J(f)$ does not contain isolated Jordan 
arcs.
\end{isoarc}
Here, by definition, a Jordan arc is called isolated (in 
$J$)
if there exists an open set which contains the arc except
for its
endpoints but no other point of $J$.
 
To prove Theorem \ref{isoarc}, we suppose that such an arc 
exists
and is parametrized by $\gamma:[0,1]\to\C$.
By Theorem \ref{rep}, the repelling periodic points are
dense in this arc. In view of Lemma \ref{la1} we may suppose
that it contains a fixed point, say, $f(z_1)=z_1$ where
$z_1=\gamma(t_1)$, $t_1\in(0,1)$.
We may also assume  that $z_1$ is not exceptional so that
$J=\overline{O^-(z_1)}$ by Lemma \ref{la4}.
Hence there exist $t_0$ and $t_2$ satisfying
$0<t_0<t_1<t_2<1$  and $n\geq 1$ such that
$f^n(\gamma(t_0))= f^n(\gamma(t_2))=z_1$
and
$f^n(\gamma(t))\not = z_1$ for $t_0<t<t_1$ and $t_1<t<t_2$.
We consider $C=f^n(\gamma([t_0,t_2]))$.
It follows from the assumption that $\gamma$ is isolated
and from the complete invariance
of $J$ that if $t\in(t_0,t_2)$,
then there exists a neighborhood $N_t$ of $f^n(\gamma(t))$
such that $N_t\cap J\backslash C=\emptyset$.
Because
$f^n(\gamma(t_0))= f^n(\gamma(t_2))=z_1$,
this is also true for $t=t_0$ and $t=t_2$.
Hence $C$ has a neighborhood $N$ satisfying
$N\cap J\backslash C=\emptyset$.
This implies that $f(C)\subset C$,
since $z_1\in C$, $J$ is completely invariant,
and $C$ is connected.
Hence $O^+(C)\subset C$. This is a contradiction,
because $O^+(C)\backslash J$ contains at most the
exceptional points of $f$ and is hence unbounded.
 
There are several other ways to see that
$J(f)=\R\cup\{\infty\}$ is impossible for $f\in E$.
One way is to observe that if this were the case,
then the upper and lower half-plane
were completely invariant with respect to $f^{\,2}$, 
contradicting
Theorem \ref{compinv1}.
Another way to see that
$J(f)=\R\cup\{\infty\}$ is impossible for $f\in E$
is to combine the complete invariance of $J(f)$
with a result of Edrei
\cite[p.\ 279]{Edr55} which says that if all roots of
$f(z)=h_n$ are real for some unbounded sequence $(h_n)$ and
an entire function $f$, then $f$ is a polynomial of degree
at most $2$.
More generally, one may use the above arguments to prove 
that
the Julia set of an entire transcendental function cannot
be contained in a finite set of straight lines; see 
\cite{Bak65}.
 
We remark that Baker, Kotus, and L\"u 
\cite[Theorem 2]{Bak91},
using a result of \v{C}e\-bo\-ta\-rev \cite{Ceb28},
have shown that if a transcendental meromorphic function $f$
satisfies
$J(f)=\R\cup\{\infty\}$, then
\[
f(z)=\varepsilon\left(cz+d+\sum_{n=1}^\infty
c_n\left(\frac{1}{a_n-z}-\frac{1}{a_n}\right)\right),
\]
where
$c,d,c_n,a_n\in \R$,
$\varepsilon=\pm 1$,
$c\geq 0$, $c_n>0$, $a_n\neq 0$,
and  $\sum_{n=1}^\infty c_n/a_n^2<\infty$.
For further details concerning
meromorphic functions satisfying $J\subset \R\cup\{\infty\}$
we refer to \cite{Bak91}.
 
\subsection{Points that tend to infinity}  \label{infinity}
Eremenko \cite{Ere89} considered the set
\[
I(f)=\{z: f^n(z)\to\infty \text{\ as\ } n\to\infty\}.
\]
If $f$ is a polynomial, then $I(f)$ is the immediate
attractive basin of the superattracting
fixed point $\infty$. In this case, we easily find that
\begin{equation}
J(f)=\partial I(f).    \label{boundaryI}
\end{equation}
Eremenko's main result in \cite{Ere89} is the following 
theorem.
\begin{ere89} \label{ere89}
If $f\in E$, then $I(f)\neq\emptyset$.
\end{ere89}
Eremenko also shows that $I(f)\cap J(f)\neq \emptyset$.
The proof of Theorem \ref{ere89} is based on the
Wiman-Valiron theory about the behavior of entire functions
near points of maximum modulus; see for example 
\cite{Hay74,Val23}.
 
Once Theorem \ref{ere89} is known, it is not difficult to
prove that (\ref{boundaryI}) holds for $f\in E$ as well.
In particular, if an entire function $f$ does not have
Baker domains (for example, if $f\in B\cap E$), then
we have $J(f)=\overline{I(f)}$.
 
We mention two questions asked by Eremenko
\cite[pp.\ 343--344]{Ere89}. Suppose $f\in E$.
\begin{qu9} \label{qu9}
Is every component of $I(f)$ unbounded?
\end{qu9}
\begin{qu10} \label{qu10}
Can every point in $I(f)$ be joined with $\infty$ by a curve
in $I(f)$?
\end{qu10}
Clearly, a positive answer to Question \ref{qu10} would 
imply that
the answer to Question \ref{qu9} is also positive.
 
Eremenko \cite[Theorem 3]{Ere89} proved that
$\overline{I(f)}$ does not have bounded components,
and he pointed out that a positive answer to Question 
\ref{qu10}
for a restricted class of functions follows from the results
of Devaney and Tangermann \cite{Dev86a}.
 
\subsection{Cantor bouquets}  \label{bouquets}
Devaney and Krych \cite{Dev84}
have studied the Julia set of exponential
functions. They find
\cite[p.\ 50]{Dev84} that if $0<\lambda<1/e$, then
$J(\lambda e^z)$ is a so-called {\em Cantor bouquet}.
We will define  Cantor bouquets below but
sketch the ideas only
briefly and refer to
\cite{Dev84a,Dev86,Dev89a,Dev91a,Dev86a} for more details.
 
For a positive integer $N$ we consider the space
$\Sigma_N$ of sequences of integers between $-N$ and $N$,
that is,
\[
\Sigma_N=\{(s_0,s_1,s_2,\dots):s_j\in\Z, |s_j|\leq N\}.
\]
There is a natural topology that makes $\Sigma_N$ into
a Cantor set. The shift $\sigma:\Sigma_N\to\Sigma_N$
is defined by 
$\sigma(s_0,s_1,s_2,\dots)=(s_1,s_2,s_3,\dots)$.
We call a closed subset $C_N$ of $\C$
a {\em Cantor-N-bouquet}
of the  meromorphic function $f$ if
$f(C_N)\subset C_N$ and if there exists
a homeomorphism
$h:\Sigma_N\times [0,\infty)\to C_N$ with the following 
properties:
\newcounter{counter2}
\begin{list}
{(\roman{counter2})}{\usecounter{counter2}}
\item $(\pi\circ h^{-1} \circ f\circ h)(s,t)=\sigma(s)$ for 
all
$t\in[0,\infty)$, where $\pi:
\Sigma_N\times [0,\infty)\to\Sigma_N$ is
the projection, that is, $\pi(s,t)=s$;
\item $\lim_{t\to\infty}h(s,t)=\infty$;
\item $\lim_{n\to\infty}f^n(h(s,t))=\infty$ if $t>0$.
\end{list}
A Cantor-$N$-bouquet is similar to a  Cantor set, but
the components are curves
tending to $\infty$ instead of points.
 
Given a sequence $C_N$ of Cantor-$N$-bouquets satisfying
$C_N\subset C_{N+1}$, the set
\[
C_\infty=\overline{\bigcup_{N=1}^\infty C_N }
\]
is called a {\em Cantor bouquet}.
 
We indicate how a Cantor bouquet can be obtained
for $E_\lambda(z)=\lambda e^z$ where $0<\lambda<1/e$.
Given $N\geq 1$, we choose $c>1$ such that
$E_\lambda(c)> c+(2N+1)\pi$ and
consider the rectangles
\[
R_j=\{z: 1<\re z<c, (2j-1)\pi < \im z<(2j+1)\pi\}
\]
for $j\in \{-N,-N+1,\dots,N\}$.
For each $j$ we have
\[
E_\lambda(R_j)=\{z: \lambda e<|z|<\lambda e^c, |\Arg z|
<\pi\}.
\]
Hence our choice of \!$c$ implies that
\!$R_k\!\subset \!E_\lambda(R_j)$ if
$j,k\!\in \!\{-N,-N+
1, \dots,N\}$.
Define $R=\bigcup_{j=-N}^N R_j$ and
\[
\Lambda_N=\{z: E_\lambda^N(z)\in R\text{\ for all\ } n\geq 1
\}.
\]
From the above observations we can deduce that $\Lambda_N$ 
is
a Cantor set homeomorphic to $\Sigma_N$.
 
This construction yields
the ``endpoints'' of the
Cantor-$N$-bouquet,
that is, the points in $h(\Sigma_N\times\{0\})$.
To obtain the curves attached to it, choose a point 
$w\in\Lambda_N$
and consider the set of all $z\in\C$ such that 
$E^n_\lambda(z)$
and $E^n_\lambda(w)$ lie in the same half-strip
\[
S_j=\{z: 1<\re z<c, (2j-1)\pi < \im z<(2j+1)\pi\}
\]
for all $n\geq 0$. This set then turns out to be a curve 
with the
desired properties. We omit the details and refer to the 
papers
cited above.
 
The method is not restricted to exponential functions. In 
fact,
it is shown \cite{Dev86a} that there exists a large class of
functions, including, for example, $\sin z$ and $\cos z$,
where the Julia set contains Cantor bouquets.
 
Besides the papers already mentioned we refer to
\cite{Aar91,Bul90,May90,Sil88} for a further discussion
of Cantor bouquets.
 
\section{Newton's method} \label{newton}
\subsection{The unrelaxed Newton method}  \label{unrelax}
Let $g$ be a meromorphic function. Newton's method
of finding the zeros of $g$ consists of iterating the
meromorphic function $f$ defined by
\begin{equation}
f(z)=z-\frac{g(z)}{g'(z)}.        \label{newtonf}
\end{equation}
In fact, if $\zeta$ is a zero of $g$, then $\zeta$ is an 
attracting
fixed point of $f$, and vice versa. The simple zeros of $g$
correspond to the superattracting fixed points of $f$.
 
Clearly, if $z$ is close enough to $\zeta$, then $f^n(z)$ 
converges
to $\zeta$ as $n\to\infty$. On the other hand, $f^n(z)$ 
cannot
tend to a zero of $g$ if $z\in J(f)$, because $J(f)$ is 
completely
invariant (under $f$). One may ask under which circumstances
it is possible that $f^n(z)$ fails to converge to zeros of
$g$ for some $z\in F(f)$ and, hence, for some open set of
$z$-values. In view of Theorem \ref{classification}
and because all fixed points of $f$ are attracting,
this is possible only in one of the following cases:
\newcounter{counter1}
\begin{list}
{(\roman{counter1})}{\usecounter{counter1}}
\item There exists $n\geq 0$ such that $f^n(z)$
is contained in a periodic cycle of immediate attractive 
basins,
Leau domains, or Siegel discs. Here the minimal
period of  the cycle is greater than $1$.
\item There exists $n\geq 0$ such that $f^n(z)$
is contained in a periodic cycle of Herman rings.
\item There exists $n\geq 0$ such that $f^n(z)$
is contained in a periodic cycle of Baker domains.
\item  $z$ is contained in a wandering domain.
\end{list}
We shall restrict here to the case that $g$ is entire and
consider the case that $g$ is a polynomial first. Then
$f$ is rational, and cases (iii) and (iv) do not occur.
It follows from a result of Shishikura \cite{Shi90}
that (ii) does not occur either. More precisely, 
Shishikura's
result says that if a rational function has only one
fixed point which is repelling or has multiplier $1$,
then its Julia set is connected.
 
On the other hand, simple examples like 
$g(z)=z^3-z+1/\sqrt{2}$
where $0$ is a superattracting periodic point of minimal 
period
$2$ for $f$ show that case (i) can occur.
From Theorem \ref{sing} and the fact that all finite fixed 
points
of $f$ are attracting we can deduce that (i) cannot occur
if $f^n(z)$ converges for all $z\in \Si(f^{-1})$.
(Theorem \ref{sing} also shows that (ii) does not occur
if $f^n(z)$ converges for all $z\in \Si(f^{-1})$.)
Because $f'(z)=g(z)g''(z)/g'(z)^2$ and $\infty$ is a fixed 
point
of $f$, we obtain the following result.
\begin{smale}  \label{smale}
Let $g$ be a polynomial, and let $f$ be defined by
\rom{(}\ref{newtonf}\rom{)}. Denote by $z_1,z_2,\dots,z_m$ 
the zeros
of $g''$ that are not zeros of $g'$. If
$f^n(z_j)$ converges for all $j\in\{1,2,\dots,m\}$, then
$f^n(z)$ converges to zeros of $g$ for all $z\in F(f)$.
\end{smale}
The proof of Theorem \ref{smale} we have sketched above
depends on Theorem \ref{sullivan}.
It is possible, however, to give a more elementary proof of
Theorem \ref{smale}.
In fact, this result can be deduced from the work of Fatou
\cite[\S 30--31]{Fat19} and Julia \cite[\S 59]{Jul18}
(see also Smale \cite[pp.\ 99--100]{Sma85}).
As an example of where
Theorem \ref{smale} applies, we mention real polynomials
with only real zeros \cite{Bar55}.
 
For a further discussion of Newton's method for polynomials 
we refer
to \cite{Hae88,Kri92,Sut92}.
If $g$ is transcendental, then so is $f$, except when 
$g=pe^q$
for polynomials $p$ and $q$. Newton's method for functions 
of
this form has been studied in detail
in~\cite{Har92}.
 
We now consider the case that $g$ and $f$ are 
transcendental.
Examples in \cite{Ber91a} show that not only (i) can occur,
but (iii) and (iv) can also occur.
In fact, it was shown that case (iii) always occurs if $g$
tends to zero in some sector sufficiently fast, for example,
if $g$ is of the form $g(z)=h(z)\exp(-z^k)$ for some 
positive
integer $k$ and some entire function $h$ of order less than 
$k$
which does not have zeros in $|\arg z|<\varepsilon$ for some
$\varepsilon >0$. An example where Newton's method leads to
wandering domains is given by
$g(z)=\exp\left(\frac{1}{2\pi i}\int_0^z\cos^2(e^u) 
\,du\right)$;
compare \cite{Ber91a}.
 
It is of interest to find classes of entire functions
for which Newton's method behaves similarly to that for
polynomials.
In \cite{Ber92a}, Newton's method for functions $g$ of the
form
\begin{equation}
g(z)=\int_0^zp(t)e^{q(t)}\,dt +c    \label{newtong}
\end{equation}
where $p$ and $q$ are polynomials and where $c$ is a 
constant
was studied. If $g$ has this form,  then $f\in N$ and hence
does not have wandering domains by Theorem \ref{nowd}.
Also,
every cycle of Baker domains contains a
singularity of $f^{-1}$ by Theorem \ref{ber92b}.
Moreover, it was shown in \cite{Ber92a}
that if $g$ and $f$ are given by
(\ref{newtong}) and (\ref{newtonf}), then $f$ does not have 
finite
asymptotic values. Hence we have the following result.
\begin{ber92} \label{ber92}
If $g$ has the form \rom{(}\ref{newtong}\rom{)} but is not 
of the form
$g(z)=e^{az+b}$ where $a$ and $b$ are constant,
then the conclusion of
Theorem \ref{smale} holds.
\end{ber92}
The case $g(z)=e^{az+b}$ has to be excluded because 
then $f(z)=z-1/a$, but we always assumed that $f$ is 
nonlinear.
In fact, the conclusion of Theorem
\ref{ber92} is false in this case.
 
Another class of entire functions where Newton's
method does not lead to wandering domains are solutions
of differential equations of the form $g''+pg=0$ where
$p$ is a polynomial. In this case, we have $f\in R$.
 
It seems likely that results of the type of Theorem 
\ref{ber92}
hold for a much wider class of functions.
For instance,
one may ask the following question which is
related to Questions \ref{qu17} and  \ref{qu5}.
\begin{qu7}
Let $g$ be a meromorphic function, 
and let $f$ be defined by
(\ref{newtonf}). Does
the convergence of $f^n(z)$
for all $z\in \Si(f^{-1})$ imply
the convergence of $f^n(z)$
(to zeros of $g$) for all $z\in F(f)$?
\end{qu7}
 
\subsection{The relaxed Newton method}   \label{relax}
As a generalization of Newton's method, one may consider the
relaxed Newton method, which is given by
iteration of
\begin{equation}
f_h(z)=z-h\frac{g(z)}{g'(z)}   \label{relaxf}
\end{equation}
where $g$ is meromorphic and $h\in\C$, $|h-1|<1$.
Again, if $\zeta$ is a zero of $g$, then $\zeta$ is an
attracting fixed point of $f$.
For $h\neq 1$, however, $\zeta$ is not
superattracting, but $f'(\zeta)=1-h/m$ if
$\zeta$ is a zero of $g$ of multiplicity $m$.
 
Clearly, the case $h=1$ corresponds to the unrelaxed Newton
method considered in the previous section. Some of the 
results
mentioned there extend to this more general case. For 
example,
if $g$ is of the form (\ref{newtong}) and if $f_h$
is defined by (\ref{relaxf}), then $f_h$ does not have 
wandering
domains, and every cycle of Baker domains of $f_h$ contains 
a
singularity of $f_h^{-1}$; see
\cite{Ber92a}.
 
The relaxed Newton method may be viewed as a discretization
of the differential equation
\begin{equation}
\dot{z}=-\frac{g(z)}{g'(z)}.    \label{flow}
\end{equation}
This differential equation has been studied in a number of
papers; compare the survey by
Jongen, Jonker, and Twilt \cite{Jon88}.
Similar to the immediate attractive basins
with respect to the iteration of $f_h$,
there are basins of attraction
with respect to the differential equation
attached to the zeros of $g$.
For a zero $\zeta$ of $g$, we
denote by $A^*(h,\zeta)$ the immediate basin of attraction
of $\zeta$ with respect to the iteration of
$f_h$, that is, the component
of $F(f_h)$ that contains $\zeta$, and
by $A(h,\zeta)$ the basin of attraction, that is,
$A(h,\zeta)=\{z:\lim_{n\to\infty}f_h^n(z)=\zeta\}$.
Clearly, $A^*(h,\zeta)\subset A(h,\zeta)$.
By $B(\zeta)$ we denote the basin corresponding to
the differential equation; that is,
$B(\zeta)$ is the set of all $w\in\C$ such that there exists
a solution
$z:[a,b)\to \C$  of (\ref{flow}) satisfying
$z(a)=w$ and $\lim_{t\to b}z(t)=\zeta$.
Considering constant solutions of (\ref{flow}), we see that
always $\zeta\in B(\zeta)$, provided we define 
$g(\zeta)/g'(\zeta)=0$
for multiple zeros $\zeta$ of~$g$.
 
We remark that if $g$ is rational, then
\begin{equation}
\Meas\left(\widehat{\C}\backslash
\bigcup_{\{\zeta:g(\zeta)=0\}}B(\zeta)\right)=0, 
 \label{measB}
\end{equation}
where $\Meas(\cdot)$ denotes Lebesgue measure on 
$\widehat{\C}$,
while
$\widehat{\C}\backslash\bigcup_{\{\zeta:g(\zeta)=0\}}
A(h,\zeta)$
may contain open sets, as already pointed out in \S 
\ref{unrelax}.
It is of interest to study
to what extent $A(h,\zeta)$
approximates $B(\zeta)$ if $h\to 0$.
For rational functions $g$ and real values of $h$
this has been studied in detail in
\cite{Hae91} and
\cite{Mei91}.
For example, it follows
from the results obtained there that
\begin{equation}
\lim_{h\to 0}\Meas\left(\widehat{\C}\backslash
\bigcup_{\{\zeta:g(\zeta)=0\}}A^*(h,\zeta)\right)=0.  
\label{measA}
\end{equation}
We also mention \cite{Fle89}, where it is proved
that if $g$ is a polynomial, then
\[
\Meas\left(
\widehat{\C}\backslash\bigcup_{\{\zeta:g(\zeta)=0\}}
A(h,\zeta)\right)=0
\]
for certain small (not necessarily real) values of $h$.
 
If $g$ is transcendental, then (\ref{measB}) and 
(\ref{measA}) need not
be true. A simple example is provided by $g=pe^q$ if $p$ and
 $q$ are
polynomials, $q$ nonconstant.
More generally, it was shown in \cite{Ber91a} that
(\ref{measB}) and (\ref{measA}) do
not hold if $g$ tends to zero in a suitable sector 
sufficiently fast.
On the other hand, we have the following result proved 
in~\cite{Ber91a}.
\begin{threlax} \label{threlax}
Let $g$ be a meromorphic function. Suppose
that $\Si(g^{-1})$ is a discrete subset of $\C$ and
that $0$ is not an asymptotic value of $g$.
Then \rom{(}\ref{measB}\rom{)} and \rom{(}\ref{measA}\rom{)}
hold.
\end{threlax}
It seems likely that the conclusion of Theorem
\ref{threlax} remains valid for more
general classes of functions.
\begin{qu8}
Is the hypothesis on the discreteness of
$\Si(g^{-1})$ necessary in
Theorem \ref{threlax}?
\end{qu8}
 
\section{Miscellaneous topics} \label{misc}
In this paper, we have concentrated on describing some
results in iteration theory that hold for all entire or
meromorphic functions or at least for large classes of
 functions.
Of course, it is  also very important
to consider
specific examples.
Already Fatou \cite[pp.\ 358--369]{Fat26} studied the 
examples
$f(z)=z+1+e^{-z}$ and $f(z)=h\sin z +a$ (where $0<h<1$ and 
$a\in \R$)
in detail, and T\"opfer \cite[\S\S 5--6]{Toe39} described 
the Julia
sets of the sine and cosine function.
 
A particularly
important topic is to consider families of functions that 
depend
on a parameter and to study how the iterative behavior 
varies as the
parameter changes.
In the iteration theory of rational functions,
the family of quadratic polynomials and its bifurcation 
diagram,
the Mandelbrot set, has been the object of much research.
Among the transcendental functions, it is probably
the exponential family
$\{\lambda e^z: \lambda\in\C\backslash \{0\}\}$
that has received most attention.
We have already mentioned some results in \S \ref{bouquets}.
 
Define $E_\lambda(z)=\lambda e^z$.
By Theorem \ref{nowd} and Corollary \ref{co3}, $E_\lambda$
does not have wandering or Baker domains.
Hence, in view of Theorem \ref{sing},
the iterative behavior of $E_\lambda$
is largely determined by the forward orbit of $0$.
In particular, $E_\lambda$ has at most one periodic cycle of
immediate attractive basins,
and if such a cycle exists, then it must contain $0$.
If $E_\lambda^n(0)\to\infty$ or if
the sequence $(E_\lambda^n(0))_{n\geq 0}$
is preperiodic but not  periodic, then 
$J(E_\lambda)=\widehat{\C}$.
 
In [25; 58; 70, \S9] 
the iteration of $E_\lambda$ has been thoroughly 
investigated;
for example, the sets
\[
D_n=\{\lambda: E_\lambda
\text{\ has an attracting periodic cycle of minimal
 period\ }
n\}
\]
have been studied in detail there.
We omit these results here but just mention
one open question.
\begin{qu11} \label{qu11}
Is $\bigcup_{n=1}^\infty D_n$ a dense subset of $\C$?
\end{qu11}
This is an analogue to a well-known conjecture of Fatou
\cite[\S 31, p.\ 73]{Fat19}
concerning rational functions.
Some partial results concerning Question \ref{qu11} can be
found in \cite{Dev84a,Dev85,Ye92,Zho89}.
 
Of course, there are many other families
that can be studied.  For example, the functions
\[
\lambda\tan z,\frac{\lambda e^z}{e^z-e^{-z}},\text{\ and\ }
\frac{e^z}{\lambda e^z+e^{-z}}
\]
were studied in \cite[\S\S 2--4]{Dev89} for certain
parameter values $\lambda$.
We also mention \cite{Dou91,How85,How90}, where numerical 
studies
concerning the iteration of transcendental meromorphic 
functions
have been carried out.
 
There are many topics that have been left out. For example,
we have not discussed ergodic problems.
There are many papers on this topic for rational maps;
see \cite[Chapter 3]{Ere90} for a survey.
Much less work has been done in this area
for transcendental functions, but we mention
\cite{Ghy85,Lyu87,Lyu88,Ree86}, which
address these questions for the exponential function.
 
Another topic we have omitted is the investigation of
the area and the Hausdorff dimension
of Julia sets of transcendental functions.
We refer to \cite[\S 7; 103; 124; 126]{Ere90a}
for results in this direction.
\section*{Acknowledgment}

I would like to thank Norbert Terglane for many useful 
discussions
and valuable suggestions. I am also grateful to Alexander 
Eremenko,
Aimo Hinkkanen, and Steffen Rohde for some helpful comments.
Finally, I am indebted  to two referees and  to the
editor, Richard Palais, for a number of helpful suggestions
that led to an improvement of this paper.

\end{document}